\documentstyle[11pt]{article}
\date{April 1996} 
\newtheorem{theo}{Theorem}  
\newtheorem{defn}[theo]{Definition}  
  
\newtheorem{lemma}{Lemma}  
\newtheorem{prop}[theo]{Proposition}  
  
\newcommand{\ssection}[1]{{\addtocounter{subsection}{1}} \vspace*{5mm}

\noindent  \thesection.\arabic{subsection} {\em  #1}.} 
  
\def\sg{{\scriptstyle {\cal G}}}
\def\End{{\mbox{\it End\/}}}

\def\mod{{\mbox{\it mod\/}}}
\def\tr{{\mbox{\it tr\/}}}
\def\be{\begin{equation}}  
\def\ba{\begin{eqnarray}}  
\def\ee{\end{equation}}  
\def\ea{\end{eqnarray}}  
\def\o{\otimes }  
\def\bo{\mbox{\,\raisebox{-0.65mm}{$\Box$} \hspace{-5.3mm}  
${\scriptstyle\times}$ \/}}  
\def\D{\Delta }  
  
\def\H{{\cal H}}  
\def\R{{\cal R}}
\def\k{{\cal k}}  
  
\def\L{{\cal L}}

\def\E{0}
\def\G{{\cal G}}  

\def\Fu{{\bf F}}  
\def\K{{\cal K}}  
\def\N{{\cal N}}  
\def\S{{\cal S}}  
  
\def\J{{\cal J}}  
  
\def\P{{\cal P}}  
\def\Z{{\cal Z}}
\def\ti{\times }  
  
\def\s{\sigma }  
\def\vs{\varsigma}

\def\a{\alpha }  
\def\b{\beta }  
\def\c{\gamma }
\def\cc{\tilde{\gamma}}  
\def\d{\delta }  
\def\e{\epsilon }

\def\k{\kappa }

\def\vth{\vartheta}

\def\vac{|0 \rangle}  
  
\def\t{\tau}  
\def\nn{\nonumber}

\newcommand{\U}[2]{\stackrel{\scriptscriptstyle #1}{U}
            \hspace*{-4mm} \phantom{U}^{#2}}
\newcommand{\T}[1]{\stackrel{\scriptscriptstyle #1}{T}
            \hspace*{-4mm} \phantom{T}^{\ }}
\newcommand{\M}[2]{\stackrel{\scriptscriptstyle #1}{M}
            \hspace*{-5mm} \phantom{M}^{#2}}
\newcommand{\Ne}[2]{\stackrel{\scriptscriptstyle #1}{N}
            \hspace*{-4mm} \phantom{N}^{#2}}
\newcommand{\Je}[2]{\stackrel{\scriptscriptstyle #1}{J}
            \hspace*{-4mm} \phantom{J}^{#2}}

\newcommand{\uE}[2]{\stackrel{\scriptscriptstyle #1}{u}
            \hspace*{-3mm} \phantom{u}^{#2}}

\newcommand{\cle}[2]{\stackrel{\scriptscriptstyle #1}{j}
            \hspace*{-2mm} \phantom{j}^{#2}}
\newcommand{\cre}[2]{\stackrel{\scriptscriptstyle #1}{j}
            \hspace*{-2mm} \phantom{j}^{#2}}
\newcommand{\cne}[2]{\stackrel{\scriptscriptstyle #1}{\eta}
            \hspace*{-2mm} \phantom{\eta}^{#2}}

\begin{document}  
\begin{titlepage}  
\title{Representation Theory of Lattice Current Algebras} 
\author{{\sc Anton Yu. Alekseev}  
\thanks{On leave of absence from Steklov Institute,  
Fontanka 27, St.Petersburg, Russia}  
\\Institute of Theoretical Physics, Uppsala University,  
\\ Box 803 S-75108, Uppsala, Sweden
\thanks{Permanent address;  e-mail:   
alekseev@rhea.teorfys.uu.se}
\\ and
\\ Institut f\"ur Theoretische Physik, 
\\ ETH -  H\"onggerberg, CH-8093 Z\"urich, Switzerland \\[4mm]  
{\sc Ludwig D. Faddeev \thanks{e-mail: faddeev@pdmi.ras.ru}}
\\ Steklov Mathematical Institute, Fontanka 27,
\\ St. Petersburg 191011, Russia \\[4mm] 
{\sc J\"urg Fr\"ohlich} 
\\ Institut f\"ur Theoretische Physik, 
\\ ETH -  H\"onggerberg, CH-8093 Z\"urich, Switzerland \\[4mm] 
{\sc Volker Schomerus \thanks{e-mail: vschomer@x4u2.desy.de}}
\\ II. Institut f\"ur Theoretische Physik, Universit\"at Hamburg,
\\ Luruper Chaussee 149, 22761 Hamburg, Germany 
\\ and
\\ Research Institute for Mathematical Sciences, Kyoto University,
\\ Kyoto 606, Japan } 
\maketitle \thispagestyle{empty}  
   
\begin{abstract}  
Lattice current algebras were introduced as 
a regularization of the left- and right moving degrees of 
freedom in the WZNW model. They provide examples of 
lattice theories with a local quantum symmetry $U_q(\sg)$. 
Their representation theory is studied in detail. In particular,
we construct all irreducible representations along with  a 
lattice analogue of the fusion product for representations 
of the lattice current algebra. It is shown that for an arbitrary 
number of lattice sites, the representation categories of the 
lattice current algebras agree with their continuum 
counterparts.  
\end{abstract}  
\end{titlepage}  
   
\section{Introduction}  
\def\tt{\tilde{\tau}}  
\setcounter{equation}{0}  
 
Lattice current algebras were introduced and first studied several 
years ago (see \cite{LKM}, \cite{FaGa} and references therein).
They were designed to provide a lattice regularization of the 
left- and right-moving degrees of freedom of the WZNW model 
\cite{WZW} and gave a new appealing view on the quantum group structure 
of the model. In spite of many similarities between lattice and 
continuum theory, fundamental relations between them remain to be 
understood. In this paper we prove the conjecture of \cite{LKM} that
the representation categories of the lattice and continuum model 
agree.   

\ssection{Lattice current algebras}
Lattice current algebras are defined over a discretized circle, i.e., 
their fundamental degrees of freedom are assigned to  $N$ vertices
and $N$ edges of a 1-dimensional periodic lattice. We enumerate 
vertices by integers $n (\mod N)$. Edges are oriented such that the 
$n^{th}$ edge points from the $(n-1)^{st}$ to the $n^{th}$ vertex. 
Being defined over lattices of size $N$, the lattice current
algebras come in families $\K_N, N$ a positive integer. A precise 
definition of these (associative *-)algebras $\K_N$ is given in 
the next section. We shall see that elements of $\K_N$ can be 
assembled into $(s \ti s)-$ matrices, $J_n$ and $N_n$, $n \in {\bf Z}
\mod N$, with $\K_N$-valued matrix elements  such that 
\ba 
  \Je{1}{}_n R \Je{2}{}_{n+1}  =  \Je{2}{}_{n+1} \Je{1}{}_n\ \ \ & , &
 \ \ \  R' \Je{1}{}_n  \Je{2}{}_n 
    =  \Je{2}{}_n \Je{1}{}_n R \ \ ,\nn \\[1mm]
  R' \Ne{1}{}_n R \Ne{2}{}_n  & = & 
 \Ne{2}{}_n R' \Ne{1}{}_n R \ \ , \label{kmex}\\[1mm]    
  \Ne{2}{}_n  \Je{1}{}_n  =  \Je{1}{}_n R \Ne{2}{}_n R' \ \ \ & , &  
  \ \ \ \Je{1}{}_n \Ne{2}{}_{n-1} = R \Ne{2}{}_{n-1} R' \Je{1}{}_n
  \ \ . \nn
\ea
To explain notations we view the matrices $J_n,N_n$ as 
elements in $\End (V) \o \K_N$, where $V$ is an s-dimensional 
vector space. In this way  $J_n = \sum m_{n,\vs}  \o l_{n,\vs}$  
determines elements $m_{n,\vs} \in \End (V)$ and $l_{n,\vs} \in \K_N$ 
which are used to define 
$$ \Je{1}{}_n = \sum m_{n,\vs} \o e \o \l_{n,\vs} \ \ , \ \ 
   \Je{2}{}_n = \sum e \o m_{n,\vs} \o \l_{n,\vs} $$ 
where $e$  is the unit in $\End (V)$. Similar definitions
apply to $N_n$. Throughout the paper we will use the symbol 
$\s$ for the permutation map $\s: \End (V)\o \End (V) \rightarrow 
\End (V) \o \End(V)$, and application of $\s$ to an object $X
\in \End (V) \o \End (V)$ is usually abbreviated by putting
a prime, i.e. $X' \equiv \s(X)$. The matrix $R = R(h) \in \End (V) 
\o \End (V)$ which appears in eqs. (\ref{kmex}) is a one-parameter 
solution of the Yang Baxter Equation (YBE). Such solutions can be 
obtained from arbitrary simple Lie algebras.       

The lattice Kac-Moody algebra $\K_N$ depends on a number of parameters, 
including the ``Planck constant'' 
$h$ in the solution $R(h)$ of the YBE and the ``lattice spacing'' $\D = 
1/N$. A first, nontrivial test for the algebraic relations (\ref{kmex})
comes from the classical continuum limit, i.e., from the limit in 
which $h$ and $\D= 1/N$ are sent to zero. Using the rules 
$$  N_n  \sim 1 - \D \eta (x)   \ \ , \ \ J_n  \sim 1 - \D j(x) \ \ , 
   \ \ R \sim 1 + i \c h r $$
with $x = n/N$, $\gamma$ being a deformation parameter
and the standard prescription 
$$    \{ . , . \} = \lim_{h \to 0} i \frac {[.,.]}{h} $$
to recover the Poisson brackets from the commutators, one finds that
\ba
\{ \cle{1}{}(x) , \cle{2}{}(y)\} &=& \frac{\c}{2} [ C, \cle{1}{}(x)-
\cle{2}{}(y) ] \d (x-y) + \gamma C \d ' (x-y) \ \ ,\nn \\[1mm]
\{ \cne{1}{}(x) , \cne{2}{}(y)\} &=& \frac{\c}{2} [ C, \cne{1}{}(x)-
\cne{2}{}(y) ] \d (x-y)  \ \ , \nn \\[1mm]
\{ \cne{1}{}(x) , \cle{2}{}(y)\} &=& \frac{\c}{2} [ C, \cle{1}{}(x)-
\cle{2}{}(y) ] \d (x-y) + \c C \d ' (x-y) \ \ . \nn 
\ea
Here $C$ is the Casimir element $C = r + r' = r + \s(r)$. For clarity, 
let us rewrite these relations in terms of components. 
When we express $C= t^a \o t^a$ and $j(x) = j_a(x) t^a$
in terms of generators $t^a $ of the classical Lie algebra, the 
relations become
\ba
\{j_a(x), j_b(y)  \} & = &\c f_{ab}^c j_c(x) \d (x-y) + \gamma \d_{ab} 
     \d ' (x-y) \nn \\[2mm]    
\{\eta_a(x), \eta_b(y)  \} & = & \c f_{ab}^c \eta_c(x) \d (x-y) \nn \\[2mm]    
\{\eta_a(x), j_b(y)  \} & = & \c f_{ab}^c j_c(x) \d (x-y) + \gamma \d_{ab} 
     \d ' (x-y) 
\nn\ea
The $f_{ab}^c\ 's$ are the structure constants of the Lie algebra, 
i.e.,  $[t^a,t^b]=f_{ab}^c t^c$.
We easily recognize the first equation as the classical Poisson 
bracket of the left currents in the WZNW model. Furthermore, the quantity 
$ j^\R(x) \equiv j^\L(x)-\eta(x)$  Poisson commutes  with $j^\L(x) 
\equiv j(x)$ 
and satisfies the Poisson commutation relations of the right 
currents, i.e.,
\ba
\{ \cre{1}{\R}(x) , \cre{2}{\R}(y)\} &=& - \frac{\c}{2} [ C, \cre{1}{\R}(x)-
\cre{2}{\R}(y) ] \d (x-y) - \gamma C \d ' (x-y) \ \ ,\nn \\[1mm]
& &\{ \cre{1}{\R}(x) , \cle{2}{\L}(y) \} = 0\ \ .  \label{Rcont}
\ea  
Hence we conclude that the lattice current algebra
as described in eqs. (\ref{kmex}) is the quantum lattice 
counterpart of the classical left and right currents.   

One would like to establish a close relationship between the 
lattice current algebra and its counterpart in the 
continuum model. A first step in this direction is described in this 
paper. We find that the representation
categories of the lattice and the continuum theory coincide. 
For this to work, it is rather 
crucial to combine left- and right-moving degrees of freedom. 
For instance, the center of the lattice current algebra
with only one chiral sector changes dramatically depending on whether
the number of lattice sites
is odd or even.  The only *-operation known for such algebras
\cite{Sklyanin-Volkov} is constructed in the case of $U_q(sl(2))$ and
does not admit straightforward generalizations.
However, no such difficulties 
appear in the full theory. It therefore appears to be rather unnatural
to constrain the discrete models to one chiral sector.  

\ssection{Main results}
In the next section we  use the methods developed in 
\cite{AGS} to provide a precise  definition of 
lattice current algebras. In contrast to the heuristic 
definition we use in this introduction, our precise formulation 
is applicable to general modular Hopf algebras $\G$, in particular
to $U_q(\sg)$, for an arbitrary  semisimple Lie algebra $\sg$. The 
main result of Section 3 provides a complete list of irreducible
representations for the lattice current algebras $\K_N$.
\\[5mm]
\noindent
{\bf Theorem A} (Representations of $\K_N$) \it 
For every semisimple modular Hopf algebra $\G$ 
and every integer $N \geq 1$, there exists a {\em lattice 
current algebra} $\K_N$ which admits a family of irreducible 
$*$-representations $D^{IJ}_N$ on Hilbert spaces $W^{IJ}_N$. 
Here the labels $I,J$ run through classes of finite-dimensional,
irreducible representations of the algebra $\G$.  \rm
\\[5mm]
The {\em two} labels $I,J$ that are needed to specify a 
representation of $\K_N$ correspond to the {\em two} 
chiralities in the theory of current algebras. In fact,
the algebra $\K_1$ is isomorphic to the quantum double
of the algebra $\G$ \cite{FN} and the pairs $I,J$ label
its representations. These results are in agreement with 
the investigation of related models in \cite{NS}.  

Next, we introduce an inductive limit $\K_{\infty}$ of the family
of finite dimensional algebras $\K_N$. It can be done using the
block-spin transformation \cite{FaGa} $$\K_N \rightarrow \K_{N+1}.$$
Under this embedding, every irreducible representation of $\K_{N+1}$ splits
into a direct sum of irreducible representations of $\K_N$. It appears that the
representation $D^{IJ}_{N+1}$ always splits into several copies of the
representation $D^{IJ}_N$.  Thus,   representations
of the inductive limit $\K_{\infty}$  are in one to one correspondence with 
representations of  $\K_1$ (or $\K_N$ for arbitrary finite $N$). 

In order to be able to take tensor products of
representations of the lattice current algebras,
we introduce a family of homomorphisms 
$$\Lambda_{N,M}:\K_{N+M-1}\rightarrow \K_N \o \K_M$$
which satisfy the co-associativity condition
$$(id\otimes \Lambda_{M,L})\circ \Lambda_{N, L+M-1}=
(\Lambda_{N,M}\otimes id)\circ \Lambda_{M+N-1,L}.$$

Let us notice that the co-product $\Lambda_{N,M}$ is supposed
to provide a lattice counterpart of the co-product defined
by the structure of superselection sectors in algebraic
field theory \cite{DHR}, \cite{DHR2}, \cite{SzVe}, \cite{Vec}.

By combining the co-product $\Lambda_{N,M}$  with the block-spin 
transformation
we construct a new co-product $\D_N: \K_N \rightarrow \K_N \otimes \K_N$
which preserves the number of sites in the lattice (see subsection 4.3).
This co-product is compatible with the block-spin transformation
and, hence, it defines a co-product for the inductive limit $\K_{\infty}$:
$$\Delta_{\infty}:\K_{\infty} \rightarrow  \K_{\infty} \otimes \K_{\infty}.$$
Our second result concerns tensor products of
representations of  the lattice current algebra $\K_{\infty}$.
\\[5mm]
\noindent
{\bf Theorem B} (Representation category of the lattice current algebra) \it 
The braided tensor categories of representations of the lattice current algebra
$\K_{\infty}$ with the co-product $\D_{\infty}$ and of the  Hopf 
algebra $\K_1$ with the co-product $\D_1$
coincide. \rm
\\[5mm]
In principle, our theory must be modified to apply to 
$U_q(\sg), q^p=1$. It is well known that 
$U_q(\sg)$ at roots of unity is not semi-simple. 
This can be cured by a process of truncation which retains
only the ``physical'' part of the representation theory 
of  quantized universal enveloping algebras. The algebraic 
implementation of this idea has been explained in \cite{AGS}
and can be transferred easily to the present situation.
We plan to propose an alternative treatment in 
a forthcoming publication.

\section{Definition of the Lattice Current Algebra}
\setcounter{equation}{0}

Our goal is to assign a family of lattice current algebras 
(parametrized by the number $N$ of lattice sites) to every 
modular Hopf-*-algebra $\G$. Before we describe the details, 
we briefly recall some fundamental ingredients from 
the theory of Hopf algebras. 
 
\ssection{Semisimple, modular Hopf-algebras}
By definition, a Hopf algebra is a quadruple $(\G,\e,\D,\S)$ of 
an associative algebra $\G$ (the ``{\em symmetry algebra} '')
with unit $e \in \G$, 
a one-dimensional representation $\e: \G \rightarrow {\bf C}$ 
(the ``{\em co-unit}''), a homomorphism $\D: \G \rightarrow \G \o \G$ 
( the ``{\em co-product}'') and an anti-automorphism $\S: \G \rightarrow 
\G$ (the ``{\em antipode}''). These objects obey a set of basic axioms
which can be found e.g. in \cite{Hopf}. The Hopf algebra $(\G,\e,\D,\S)$ 
is called quasitriangular if there is an invertible element $R \in 
\G \o \G$ such that  
\ba
     R\   \D (\xi) & = & \D'(\xi) \ R \ 
      \ \mbox{ for all } \ \ \xi \in \G\ , \nn \\[1mm]
     (id \o \D)(R) = R_{13}R_{12} \ \ & , & \ \ 
     (\D \o id)(R) = R_{13}R_{23} \ \ .  \nn 
\ea
Here $\D' = \s \circ \D$, with $\s: \G \o \G \rightarrow \G \o \G$ the 
permutation map, and we are using the standard notation for the elements 
$R_{ij} \in \G \o \G \o \G$. 
 
For a ribbon Hopf-algebra one postulates, in addition, the existence of 
a certain invertible central element $v \in \G$ (the ``{\em ribbon element}'')
which factorizes $R'R \in \G \o \G$ ( here $R' = \s(R)$), in the sense that 
$$  R'R = (v \o v) \D(v^{-1})   \ \  $$ 
(see \cite{ReTu} for details). The ribbon element $v$ and the 
element $R$ allow us to construct a distinguished grouplike 
element $g \in \G$ by the formula 
   $$  g^{-1} = v^{-1} \sum \S(r^2_\vs) r^1_\vs \ 
       \ \ , $$
where the elements $r^i_\vs$ come from the expansion 
$R = \sum r^1_\vs \o r^2_\vs$
of  $R$. The element $g$ is important in the definition 
of q-traces below. 

We want this structure to be consistent with a *-operation
on $\G$. To be more precise, we require that 
\be 
R^* = (R^{-1})' = \s(R^{-1}) \ \  ,\ \ \D(\xi)^* = \D'(\xi^*),
\label{stop}
\ee
and that $v,g$ are unitary \footnote{Here we have fixed 
$*$ on $\G \o \G$ by $(\xi \o \eta)^* = \xi^* \o \eta^*$. 
Following \cite{MS}, we could define an alternative involution 
$\dagger$ on $\G \o \G$ which incorporates a permutation of 
components, i.e., $(\xi \o \eta)^\dagger = \eta^\dagger \o 
\xi^\dagger$ and $\xi^\dagger = \xi^* $ for all $\xi, \eta \in \G$. 
With respect to $\dagger$, $\D$ becomes an ordinary $*$
-homomorphism and $R$ is unitary.}.  
This structure is of particular interest, since it appears in the theory of 
the quantized universal enveloping algebras $U_q(\sg)$ when the complex
parameter $q$ has values on the unit circle \cite{MS}. 

At this point we assume that {\em $\G$ is 
semisimple}, so that every representation of $\G$ can
be decomposed into a  direct sum  of finite-dimensional, irreducible
representations. From 
every equivalence class $[I]$ of irreducible representations of $\G$, 
we may pick a representative $\t^I$, i.e., an irreducible representation
of $\G$ on  a $\d_I$-dimensional Hilbert space $V^I$. The {\em quantum
trace} $tr_q^I$ is a  linear functional acting on elements $X \in 
\End (V^I)$ by 
$$     tr^I_q (X) = Tr^I( X \t^I(g))\ \ . $$
Here $Tr^I$ denotes  the standard trace on $\End (V^I)$ with 
$Tr^I (e^I) = \d_I$ and $g \in \G$ has been defined above. Evaluation
of the unit element $e^I \in \End (V^I)$ with $tr^I_q$ gives the 
{\em quantum dimension} of the representation $\t^I$, 
$$    d_I  \equiv tr^I_q(e^I)     \ \  . $$
Furthermore, we  assign a  number $S_{IJ}$ to every pair of 
representations $\t^I,\t^J$, 
$$    S_{IJ} \equiv  \N (tr^I_q \o tr^J_q) (R'R)^{IJ} 
  \ \ \ \mbox{ with } \ \ (R'R)^{IJ} =(\t^I \o \t^J)(R'R) \ \ , $$
for a suitable, real normalization factor $\N$. 
The numbers $S_{IJ}$ form the so-called {\em S-matrix} $S$. 
Modular Hopf algebras are ribbon 
Hopf algebras with an  invertible $S$-matrix \footnote{If
a diagonal matrix $T$ is introduced according to $T_{IJ} = 
\varpi \d_{I,J} d_I^2 \t^I(v)$ (with an appropriate 
complex factor $\varpi$), then $S$ and $T$ furnish a 
projective representation of the modular group $SL(2,{\bf Z})$.}. 

Let us finally recall that the  tensor 
product, $\t \bo \t'$, of two representations $\t,\t'$ of a Hopf 
algebra is defined by 
$$ (\t \bo  \t')(\xi) = (\t \o \t') \D (\xi) \ \ \mbox { for all } 
 \ \ \xi \in \G\ \ . $$
In particular, one may construct the tensor product $\t^I \bo \t^J$ 
of  two irreducible representations. According to  our assumption  
that $\G$ is semisimple, such tensor products of representations can 
be decomposed 
into a direct sum  of irreducible representations, $\t^K$. The 
multiplicities
$N^{IJ}_K$ in this Clebsch-Gordan decomposition of $\t^I \bo \t^J$ 
are called {\em fusion  rules}.

Among all our assumptions on the structure of the Hopf-algebra $(\G,
\e,\D,\S)$
(quasi-triangularity, existence  of a ribbon  element $v$, 
semisimplicity of $\G$ and invertibility of $S$), semisimplicity of 
$\G$ is the most problematic one. In fact it is violated by the 
algebras $U_q(\sg)$ when $q$ is  a root of unity. It is 
sketched in \cite{AGS} how ``truncation'' can cure this problem, 
once the theory has been extended to weak quasi-Hopf algebras
\cite{MS}. 
\vspace*{5mm}

\noindent
{\sc Example:} {\em (Hopf-algebra $\Z_q$)} 
We wish to give one fairly trivial example for the 
algebraic structure discussed so far. Our example comes from 
the  group ${\bf Z}_p$. To be more precise, we consider
the  associative algebra $\G$ generated by one element $g$ subject 
to the relation $g^p = 1$. On this algebra, a co-product, 
co-unit and an antipode  can be defined by 
$$ \D(g) = g \o g \ \ , \ \ \S (g) =  g^{-1} \ \ , \ \ 
 \e(g) = 1 \ \  . $$
We observe that  $\G$ is a commutative semisimple algebra. It has 
$p$ one-dimensional representations $\t^r(g) = q^r, r = 0, \dots, p-1,$
where $q$ is a root of unity, $q = e^{2 \pi i/p}$. We may 
construct characteristic projectors $P^r \in \G$ for these 
representations according to 
$$ P^r = \frac{1}{p} \sum_{s=0}^{p-1}  q^{-rs}  g^s\ \ \mbox{ for }
\ \ r= 0, \dots, p-1\ \  . $$
One can easily check that $\t^r(P^s) = \delta_{r,s} $.  
The elements $P^r$ are employed to obtain a nontrivial $R$-matrix,  
$$ R = \sum_{r,s} q^{rs} P^r \o P^s \ \ . $$ 
When evaluated with a pair of representations $\t^r, \t^s$ we find
that $(\t^r \o \t^s)(R) = q^{rs}$. The $R$-matrix satisfies all 
the  axioms stated above and thus turns $\G$ into a quasitriangular 
Hopf-algebra. Moreover, a ribbon element is provided by  $v = 
\sum q^{-r^2} P_r$. We can finally introduce a $*$-operation on $\G$ 
such that $g^*= g^{-1}$. The consistency relations \ref{stop} 
follow from the co-commutativity of $\D$, i.e. $\D' = \D$, and
the property $R = R'$. A direct computation 
shows that the $S$-matrix is invertible only for  odd integers $p$. 
Summarizing all this, we have constructed a family of semisimple  
ribbon Hopf-*-algebras $\Z_q, q= exp(2\pi i /p) $. They are 
modular Hopf-algebras for all odd integers $p$.  

\ssection{$R$-matrix formalism}
Before we propose a definition of lattice current algebras, we  
mention that Hopf algebras $\G$ are intimately related to 
the objects $N_n$, $n \in {\bf Z}\mod N,$ introduced in eq. (\ref{kmex}).
To understand this relation, let us introduce another (auxiliary)
copy, $\G_a$, of $\G$ and let us consider the $R$-matrix as an   
object in $\G_a \o \G$. To distinguish the latter clearly from 
the usual $R$, we denote it by $N_\pm$,  
$$ N_- \equiv R^{-1} \in \G_a \o G \ \ \ , \ \ \ N_+ \equiv R' 
\in \G_a \o \G\ \ . $$
Quasi-triangularity   
of the $R$-matrix furnishes the relations
\ba \label{Npmeq}
 \D_a(N_{\pm}) =  \Ne{2}{}_{\pm} \Ne{1}{}_{\pm} \ \ & , & \ \  
 R \Ne{2}{}_+  \Ne{1}{}_- = \Ne{1}{}_- \Ne{2}{}_+ R\ \ , \\[1mm]
 R \Ne{2}{}_\pm \Ne{1}{}_\pm  & = & 
   \Ne{1}{}_\pm \Ne{2}{}_\pm  R \ \ . \nn 
\ea 
Here we use the same notations as in the introduction, and 
$\D_a (N_\pm) = (\D \o id) (N_\pm) \in \G_a \o \G_a \o \G$. 
The subscript $\ _a$ reminds us that $\D_a$ acts on the 
auxiliary (i.e. first) component of $N_\pm$. To be perfectly 
consistent, the objects $R$ in the preceding equations should 
all be equipped with a lower index $a$ to show that $R \in \G_a 
\o \G_a$ etc.. We hope that no confusion will arise from omitting 
this subscript on $R$. The equations (\ref{Npmeq}) are somewhat 
redundant: in fact, the exchange relations on the second line 
follow from the first equation in the first line. This 
underlines that the formula for $\D_a(N_\pm)$ encodes information 
about the product in $\G$ rather than the co-product
\footnote{The co-product $\D$ of $\G$ acts on $N_\pm$ according 
to $\D(N_\pm)= (id \o \D)(N_\pm) = N_\pm \tilde N_\pm \in 
\G_a \o \G \o \G$. Here $N_\pm \ [\tilde N_\pm]$ on the right 
hand side of the equation have the unit element $e \in \G$ in 
the third \ [second] tensor factor}. More explanations of this 
point follow in Subsection 2.3.  

Next, we combine $N_+$ and $N_-$ into one element  
$$   N \equiv N_+ (N_-)^{-1} = R'R  \ \ \in \ \G_a \o \G\ \ .$$
{}From the properties of $N_\pm$ we obtain an expression for the 
action of $\D_a$ on $N$,  
\ba 
  \D_a(N)  & = & \Ne{2}{}_+ \Ne{1}{}_+ (\Ne{1}{}_-)^{-1} 
             (\Ne{2}{}_-)^{-1}\nn \\[1mm]
         & = & R^{-1} \Ne{1}{}_+ \Ne{2}{}_+ R (\Ne{1}{}_-)^{-1} 
             (\Ne{2}{}_-)^{-1} \nn\\[1mm]
         & = & R^{-1} \Ne{1}{}_+ (\Ne{1}{}_-)^{-1} R \Ne{2}{}_+ 
             (\Ne{2}{}_-)^{-1} \nn\\[1mm]
         & = & R^{-1} \Ne{1}{} R \Ne{2}{}\ \ . \nn
\ea
As seen above, the formula for $\D_a(N)$ encodes relations in the
algebra $\G$ and  implies, in particular, the following exchange 
relations for $N$:
\ba
    R' \Ne{1}{} R \Ne{2}{} 
    &=&  R'R \D_a (N)  = R' \D'_a(N) R \nn \\[1mm]
    &=& \Ne{2}{} R' \Ne{1}{} R \ \ . \nn 
\ea
This kind of  relations first appeared in \cite{ReSTS}.
We used them in the introduction when 
describing the objects $N_n $ assigned to the sites of the 
lattice. Thus we have shown that, for any modular 
Hopf-algebra $\G$, one may construct objects $N$ obeying the  
desired quadratic relations.

The other direction, namely the problem of how to construct a modular 
Hopf-algebra $\G $ from an object $N $ satisfying the exchange 
relations described above, is more subtle.  
To begin with, one has to choose linear
maps $\pi: \G_a \rightarrow {\bf C}$ in the dual $\G'_a$ of $\G_a$. When     
such linear forms $\pi \in \G'_a$ act on the first tensor factor
of $N \in \G_a \o \G$ they produce elements in $\G$:         
$$        
          \pi(N)  \equiv (\pi\o id) (N) \ \ \in \ \G \ \ 
          \mbox{ for all } \ \ \pi \in \G'_a\ \ . 
$$ 
$\pi(N) \in \G$ will be called the $\pi$-component of $N$ or just 
{\em component of $N$}. Under certain technical assumptions it has been
shown in \cite{AlSc} that the components of $N$ generate the 
algebra $\G$. In this sense one can reconstruct the modular 
Hopf-algebra $\G$ from the object $N$. 

\begin{lemma} \cite{AlSc} Let $\G_a$ be a finite-dimensional, 
semisimple modular Hopf algebra and $\N$ be the algebra generated 
by components of $N \in \G_a \o \N$ subject to the relations 
$$  \Ne{1}{} R \Ne{2}{} = R \D_a(N) \ \ ,  $$ 
where we use the same notations as above. Then $N$ can be 
decomposed into a product of elements $N_\pm \in \G_a \o \N$, 
\ba
 N  =  N_+ N_-^{-1}  & & \mbox{such that} \nn \\[3mm]
 \D(N) \equiv  N_+ \tilde N N_-^{-1} & & \in \G_a \o \N \o \N  
 \nn \\[1mm]
 \e(N) \equiv  e \in \G_a \ & , & \  \S(N_\pm)  \equiv  
 N_\pm^{-1} \in \G_a \o \N\nn 
\ea
define a Hopf-algebra structure on $\N$. Here, the action of 
$\D,\e,\S$ on the second tensor component of $N,N_\pm$ is 
understood. In the equation for $\D(N)$, $N_+,N_-^{-1}$ are
supposed to have a trivial entry in the third component while
$\tilde N = \sum n_\vs \o e \o N_\vs $ with $e$ being the unit 
in $\N$ and $N = \sum n_\vs \o N_\vs \in \G_a \o \N$.     
As a Hopf algebra, {\em $\N$ is isomorphic to $\G_a$}. 
\end{lemma}

Let us remark that the $*$-operation in $\G$ induces a $*$-operation
in $\N$ which looks as follows:
\be
N_+^*=N_-.
\ee

In our definition of the lattice current algebras below, 
we shall describe the degrees of freedom at the lattice sites 
directly in terms of elements $\xi  \in \G$, instead of working 
with $N$ (as  in the introduction).  

The $\d_I$-dimensional representations $\t^I$ of $\G_a \cong \G$ furnish 
a $\d_I \ti \d_I$-matrix of linear forms on $\G_a$. When these forms 
act on the first tensor factor of  $N$, we obtain  a matrix 
$N^I \in \End(V^I) \o \G$ of  elements in $\G$, 
$$    N^I \equiv \t^I (N) = (\t^I \o id)(N)\ \ . $$
These matrices will turn out to be useful.  
\vspace*{5mm}

\noindent
{\sc Example:} {\em ($R$-matrix formalism for $\Z_q$)}
Let us illustrate all these remarks on  the example of 
$\G = \Z_p$. Recall that $R = \sum q^{rs} P^r \o P^s$ and 
that $\Z_p$ has only one-dimensional representations 
given by $\t^r(g)  = q^r$. Evaluation of the 
objects $N_{\pm}$ in representations $\t^r$ produces 
elements  $N^r_{\pm} = (\t^r \o id)(N_{\pm}) \in 
{\bf C } \o \Z_q \cong \Z_q$.  Explicitly, they 
are given by 
$$ N^r_+ = \sum q^{rs} P_s = g^r  \ \ \mbox{ and } \  \ 
   N^r_- = \sum q^{-rs} P_s = g^{-r} \ \ .$$
Together with the property $(\t^r \o \t^s)\D \cong \t^{r+s}$ 
the  relations (\ref{Npmeq}) become
$$  g^{\pm(r+s)} = g^{\pm s} g^{\pm r} \ \ ,  \  \  
    q^{rs} g^s g^{-r} = g^{-r} g^s q^{rs}\ \ . $$
For $N \in \Z_q \o \Z_q$ we find
$$     N  = \sum q^{2rs}P_r \o P_s \ \  \mbox{ and } \ \ 
       N^r = g^{2r}\ \ .  $$
As predicted by the general theory, the elements 
$N^r \in \Z_q$ generate the  algebra $\Z_q$ when $p$ is odd.  

\ssection{Definition of $\K_N$} 
Next, we turn to the definition  of 
the lattice current algebras $\K_N$ associated to a fixed 
modular Hopf algebra. Before entering the abstract formalism, 
it is useful to analyse the classical geometry of the discrete model. 
Our classical
continuum theory contains two Lie-algebra valued fields, namely 
$\eta(x)$ and $j(x)$. To describe a configuration of $\eta$, 
for instance, we have to place a copy of the Lie-algebra at 
every point $x$ on the circle. On the lattice, there are only
$N$ discrete points left and hence configurations of the lattice
field $\eta$ involve only $N$ copies of the Lie algebra. When 
passing from the continuum to the lattice, we encode the 
information about the field $j(x)$ in the holonomies along 
links, 
$$  j_n =P exp(\int_n j(x) dx) \ \ . $$ 
Here $\int_n$ denotes integration along the $n^{th}$ link that
connects the $(n-1)^{st}$ with the $n^{th}$ site. 
The classical lattice field $j_n$ has values in the Lie group. 
Let us remark that, even at the level of Poisson brackets, 
the variables $j_n$ can not be easily included into the Poisson algebra.
The reason is that $j_n$'s fail to be continuous functions of the
currents. Therefore, we should regularize the Poisson brackets
(or commutation relations) of the lattice currents. This regularization
is done in the most elegant way with the help of  $R$-matrices.
This consideration explains an immediate appearance of the
quantum groups in the description of the lattice current algebras.
 
In analogy to the classical description of the lattice field $\eta$, 
the lattice current algebras contain $N$ commuting copies of 
the algebra $\G$ or, more precisely, $\K_N$ contains an $N$-fold 
tensor product $\G^{\o_N}$ of $\G$ as a subalgebra. We denote 
by $\G_n$ the  subalgebra 
$$ \G_n = e \o \dots \o \G \o \dots \o e \ \subset \G^{\o_N}$$
where $\G $ appears in the $n^{th}$ position and all other entries 
in the tensor product are trivial. The canonical isomorphism
of $\G$ and $\G_n \subset \G^{\o_N}$ furnishes the homomorphisms
$$ \iota_n:\G \rightarrow \G^{\o_N} \ \ \mbox{for all} \ \ 
   n = 1, \dots, N .  $$ 
We think of the copies $\G_n$ of $\G$ as being placed at the 
$N$ sites of a periodic lattice, with $\G_n$ assigned to the 
$n^{th}$ site. In addition, the definition of $\K_N$ will 
involve generators $J_n, n= 1,\dots,N$. The generator $J_n$ 
sits on the link connecting the $(n-1)^{st}$ with the $n^{th}$ site.

\begin{defn} The lattice current algebra $\K_N$ is generated
by components \footnote{Recall that a component of $J_n$ is an element 
$\pi(J_n) \equiv (\pi \o id) (J_n)$ in the algebra $\K_N$. Here   
$\pi$ runs through the dual $\G'_a$ of $\G_a$.} 
of $J_n \in \G_a \o \K_N, n = 1, \dots, N, $ 
along with elements in $ \G^{\o_N}$. These generators
are subject to three different types of relations.  
\begin{enumerate}
\item  {\em Covariance} properties express that the $J_n$ are tensor 
       operators transforming under the action of 
       elements $\xi_m \in \G_m$ like holonomies in a 
       gauge theory, i.e., 
        \ba
    \iota_n(\xi) J_n &=& J_n \D_n (\xi) \ \ \mbox{ for all }
    \ \ \xi \in \G \nn\\
    \D_{n-1}( \xi ) J_n &=& J_n \iota_{n-1}(\xi) \ \ \mbox{ for all }
    \ \ \xi \in \G \label{cov}\\
    \iota_m(\xi) J_n&=& J_n \iota_m(\xi) \ \ \mbox{ for all }
    \ \ \xi \in \G , m \neq n,n-1\  \mod N\ \ . \nn
   \ea
  The covariance relations (\ref{cov}) make sense as relations
  in $\G_a \o \K_N$, if $\iota_n(\xi) \in \G_n \subset \K_N$ is 
  regarded as an element $\iota_n(\xi) \in \G_a \o \K_N$ with 
  trivial entry in the first tensor factor and $\D_n(\xi) 
  \equiv (id \o \iota_n ) \D (\xi) \in \G_a \o \G_n \subset 
  \G_a \o \K_N$.
\item {\em Functoriality } for elements $J_n$ on a fixed
  link means that 
  \be
   \label{OPE}
        \Je{1}{}_n \Je{2}{}_n =  R \D_a( J_n) 
 \ee
  This is to be understood as a relation in $\G_a \o \G_a \o \K_N$ 
  where $\D_a : \G_a \o \K_N \rightarrow \G_a \o \G_a \o \K_N$ acts trivially
  on the  second tensor factor $\K_N$ and $R = R \o e \in \G_a 
  \o \G_a \o   \K_N$. The other notations were explained in the 
  introduction. We also require that the elements $J_n$ possess 
  an inverse   $J_n^{-1}   \in \G_a \o \K_N$  such that 
  \be
   \label{invers}
   J_n J_n^{-1}  =  e  ,  \ \ \
   J_n^{-1} J_n = e\ \ .
  \ee
\item
  {\em Braid relations} between elements $J_n, J_m$ assigned
  to different links have to respect the gauge symmetry and locality
  of the model. These principles require
  \ba
  \Je{1}{}_n \Je{2}{}_m =  \Je{2}{}_m \Je{1}{}_n & & 
  \ \mbox{ if } n \neq m,m \pm 1 \mod N, \nn \\
  \Je{1}{}_n R \Je{2}{}_{n+1} & =&  \Je{2}{}_{n+1} \Je{1}{}_n
  \label{braid}
  \ea
  $R$ denotes the element $R \o e \in \G_a \o \G_a \o  \K_N$ as before.
\end{enumerate}
\end{defn}

The lattice current algebra $\K_N$ contains a subalgebra $\J_N$
generated by components of the $J_n$ only. They are subject to 
functoriality (2.) and braid relations (3.). The subalgebra $\J_N$
admits an action of $\G^{\o_N}$ (by generalized derivations)
such that the full lattice current algebra $\K_N$ can be regarded 
as a semi-direct product of $\J_N$ and $\G^{\o_N}$ with respect 
to this action. Our covariance relations (1.) give a precise 
definition of the semi-direct product. 
 
Let us briefly explain how  Definition 1 is related  to the 
description we used in the introduction. The relation  of the 
Hopf algebras $\G_n$ and the objects $N_n$ has already 
been discussed. Our covariance relations in 
eq. (\ref{cov}) correspond to the exchange relations between 
$N$ and $J$ in the third line of eq. (\ref{kmex}). They can 
be related explicitly with the help of the  quasi-triangularity 
of $R$, using the formula $N = R'R$. We have, for instance,  
\ba
   \Ne{2}{}_n \Je{1}{}_n &=& (e \o (R'R)_n) \Je{1}{}_n 
    = \Je{1}{}_n \ [(id \o \D_n)(R'R)]_{213} \nn \\[1mm]  
   & =&   \Je{1}{}_n\  R (e \o (R'R)_n) R' = \Je{1}{}_n R \Ne{2}{}_n R'
   \ \ , \nn
\ea
where we use the notation $(R'R)_n = (id \o \iota_n)(R'R) 
\in \G_a \o \G_n $ and $R = (R \o e) \in \G_a \o \G_a \o \K_N$ as  
usual. $[. ]_{213}$ means that the first and the second tensor factors
of the expression inside the brackets are exchanged. For 
finite-dimensional semisimple modular Hopf algebras $\G$, 
Lemma 1 implies that one could define $\K_N$ using the objects 
$N_n \in \G_a \o \G_n \subset \G_a \o \K_N$ instead of elements 
$\eta \in \G^{\o_N}$. The generators $N_n$ would have to obey 
the exchange relations stated in eq. (\ref{kmex}) and 
$$ \Ne{1}{}_n R \Ne{2}{}_n = R \D_a(N_n)\ \ . $$

The functoriality relation for $J_n$ did not appear in the 
introduction. But we can use it now to derive quadratic relations
for the $J_n$ in much the same way as  has been done for $N$, 
earlier in this section. Indeed we find
\ba
  R' \Je{1}{}_n \Je{2}{}_n & = & R'R \D_a (J_n) = R' \D_a'(J_n) 
   R \nn \\[1mm]
   & = & \Je{2}{}_n \Je{1}{}_n \ R \nn \ \ . 
\ea
This exchange relation is the one used in the introduction to describe 
the lattice currents $J_n$. When the  first two  tensor factors 
in this equation are evaluated with representations of $\G_a \cong \G$ 
one derives quadratic relations for the $\K_N$-valued matrices 
$$  J_n^I \equiv (\t^I \o id)(J_n) \in \End(V^I) \o \K_N \ . $$

As we discussed earlier in this  section, elements in the algebra
$\K_N$ can be obtained from $J_n$ with the help of  linear forms 
$\pi \in \G'_a$. To understand {\em functoriality} properly one must 
realize that it describes the  ``multiplication table'' for  
elements $\pi(J_n)  \in \K_N$. If we pick two linear forms 
$\pi_1, \pi_2  \in \G'_a$ on $\G_a$,  the  corresponding elements
in $\K_N$ satisfy
$$  \pi_1(J_n) \pi_2(J_n) = (\pi_1 \o  \pi_2)(R\D(J_n))\  \  .$$
We can rewrite this equation by means of the ({\em twisted}) 
associative product $\ast$ for elements $\pi_i \in \G'_a$, 
$$  (\pi_1 \ast \pi_2) (\xi) \equiv (\pi_1 \o \pi_2)
    (R\D(\xi))\ \  . $$
It allows us to express the product $\pi_1(J_n) \pi_2(J_n)$ 
in terms of the element $\pi_1 \ast \pi_2 \in \G_a'$, 
$$ \pi_1(J_n)  \pi_2(J_n) = (\pi_1  \ast \pi_2)(J_n)\  \   
\mbox{ for all } \ \ \pi_1, \pi_2 \in \G'_a\ \ . $$

\noindent {\bf Remark:} It may help here to invoke the analogy with 
a more familiar situation. In fact, the Hopf algebra structure of 
$\G$ induces the standard (non-twisted) product  $\cdot$ 
on its dual $\G'$,  
$$ (\pi_1 \cdot \pi_2) (\xi) \equiv (\pi_1 \o \pi_2) (\D(\xi))
   \  \ \mbox{  for all } \ \ \xi \in \G \ \ . $$
Let us identify $\pi\in\G'_a$ with the image $\pi(T)=(\pi\o id)(T)$ 
of some universal object $T \in \G_a \o \G'$ and insert $T$ 
into the definition of the product $\cdot$,  
$$  \pi_1 (T) \pi_2  (T) \equiv (\pi_1 \o \pi_2 )(\D_a( T))
      = (\pi_1 \cdot \pi_2)(T)  \ \  . $$
Here and in the following we omit the $\cdot$ when multiplying 
elements $\pi(T)$. The derived multiplication rules for components 
of $T \in \G_a  \o \G'$  are equivalent to the following functoriality 
 \ba  \T{1} \T{2} &  =  & \D_a (T) \ \nn \\[1mm] 
 \mbox{and imply {\it RTT-relations}:} & & 
  R \T{1} \T{2}  =   \T{2}  \T{1} R \ \ \  .  \nn 
\ea 
Such relations are known to define quantum groups, and hence our 
variables $J_n$ describe some sort of twisted quantum groups. This 
fits nicely with the nature of the classical lattice field $j_n$.
As we have noted earlier, the latter takes values in a Lie group. 

We saw above that the definition of a product for components 
of $J_n$ implies the desired quadratic relations. The converse 
is not true in general, i.e. functoriality is a stronger 
requirement  than the quadratic relations. In the familiar case
of $U_q(sl_2)$ for example, functoriality furnishes also the 
standard determinant relations which are usually ``added by hand''
when algebras are defined in terms of quadratic relations. Due 
to functoriality we are thus able to develop a universal theory 
which does not explicitly depend on the specific properties of 
the Hopf algebra $\G$.

We have shown that the mathematical definition of $\K_N$ 
represented in this section agrees with 
the description used in the introduction. The algebra 
$\K_N$ now appears as a special example of the ``graph algebras'' 
defined and studied in \cite{AGS} to quantize Chern-Simons theories. 
This observation will enable us to use some of the general properties 
established there.             

The most important one among such general properties is the existence 
of a *-operation. On the copies $\G_n$, a *-operation comes from 
the structure of the modular Hopf-*-algebra $\G$. Its action 
can be extended to the algebra $\K_N$ by the formula
\begin{equation} \label{inv}
J_n^* = S_n^{-1} J_n^{-1} S_{n-1} \ \ , 
\end{equation}
where $S_n =(id \o \iota_n)[\D(\k^{-1}) (\k \o \k) R^{-1}] \in 
\G_a \o \G_n \subset \G_a \o \K_N$. Here $\iota_n : \G \rightarrow \G_n 
\subset \K_N$ is the canonical embedding, and the central element 
$\k \in \G$ is  a certain square root of the  ribbon  element 
$v$, i.e. $\kappa^2 = v$ (cp. \cite{AGS} for details). Let us note that the
formula (\ref{inv}) can be rewritten using  the elements
$N_{n,\pm} \in \G_a \o \G_n$ constructed from the $R$-element
according to 
$$      N_{n,+} \equiv (id \o \iota_n)(R') \ \ , \ \ 
     N_{n,-} \equiv (id \o \iota_n)(R^{-1})   \ \ . $$
The conjugate current is expressed as
\begin{equation}
J_n^*=(\kappa_{n-1} \kappa_n)^{-1} \ N_{n,-}^{-1} J_n^{-1} N_{n-1,-} \
(\kappa_{n-1} \kappa_n) \ \ .
\end{equation}
Here $\kappa_{n-1}=\iota_{n-1}(\kappa)$, $\kappa_{n}=\iota_{n}(\kappa)$.
In order to verify the property $(J_n^*)^*=J_n$, one uses the following
identities:
$$
v_{n-1}^{-1} J_n v_{n-1}= v_a N_{n-1}^{-1} J_n \ , \  v_n^{-1} J_n v_n=v_a^{-1} 
J_n N_n,
$$
where $v_a$ is the ribbon element in the auxillary Hopf algebra $\G_a$.
The ribbon elements $v$ at different lattice sites generate automorphisms of the
lattice current algebra which resemble the evolution automorphism in the quantum
top \cite{A-F}.

\vspace*{5mm}

\noindent
{\sc  Example:}{\em (The $U(1)$-current algebra on the  lattice)}
It is  quite instructive to apply the general definition  of 
lattice current algebras to the  case $\G = \Z_q$. 
Recall that $\Z_q$ is generated by one unitary element $g$ 
which satisfies $g^p =1$. Representations $\t^s$ of $\G$ were 
labeled by an integer $s = 1, \dots, p-1,$ and 
$\t^s(g) = q^s$ with $q= \exp(2\pi i/ p)$. We can
apply the one-dimensional representations $\t^s$ to the 
current $J_n \in \G_a \o \K_N$ to obtain elements $J^s_n =\t^s(J_n) = 
(\t^s \o id)(J_n) \in \K_N$. Functoriality becomes 
$$      J^s_n J^t_n  = q^{ts} J^{s+t}_n \ \ , $$
where we have used that $(\t^s \o \t^t)\D_a (\xi) = \t^{s+t}(\xi)$ 
and $(\t^s \o \t^t)(R^{-1}) = q^{-st}$ ($s+t$ is to be understood 
modulo $p$). The relation allows to generate the elements $J^s_n$
from $J^1_n \in \K_N$. Observe that $J^0_n$ is 
the unit  element $e$ in the algebra $\K_N$. From the previous 
equation  we deduce that the $p^{th}$ power of the generator 
$J^1_n$ is proportional to $e$, 
$$  (J^1_n)^p =  q^{p(p-1)/2} J^p_n =  q^{p(p-1)/2} e\ \  .$$  
This  motivates us to  introduce the renormalized generators 
$w_n = q^{(1-p)/2} J^1_n \in \K_n$ which  obey $w_n^p  =  e$.  
In the following it suffices to specify relations for the generators 
$w_n$ and $g_n = \iota_n(g)$ of $\K_N$. The covariance equations
(\ref{cov}) read
$$    g_n w_n = q w_n g_n \ \ , \ \ w_n g_{n-1} = q g_{n-1} w_n$$
since $(\t^1 \o id)\D(g) = q g$. The exchange relations for    
currents become 
$$     w_n w_{n+1} = q^{-1} w_{n+1} w_n \ \ . $$
The identity $\D(\k^{-1}) (\k \o \k)R^{-1}= e \o e$ with $\k  = \sum 
q^{-r^2/2}P_r$ finally furnishes 
$$    w_n^* = w_n^{-1}\  \  .$$     
At this point one can easily recognize the algebra of $w_n$'s  
as the lattice $U(1)$-current algebra \cite{FaVo}. 

\ssection{The right currents}
Let us stress that there is a major ideological  difference 
between our discussion of lattice current algebras and the 
work in \cite{AGS}. In 
the context of Chern Simons theories, the graph algebras 
were introduced as  auxiliary objects, and physical variables 
of the theory were to be constructed from objects assigned to 
the links, i.e. from the variables $J_n$. Here the $J_n$'s 
represent only the left-currents, and we expect also right- 
currents to be present in the theory. They are constructed 
from the variables  $J_n$ {\em and } the elements $\xi_n \in \G_n$ 
and hence give a physical meaning to the graph algebras.  
We define a family of new variables $J^\R_n \in \G_a \o \K_N$ 
on the lattice by setting
$$ J^{\R}_n =  N_{n,-}^{-1} J_n^{-1} N_{n-1,+}
    \ \ . $$
The $J^{\R}_n$ turn out to  provide the right 
currents in our theory. For the rest of this section we will 
use the symbol $J^{\L}_n$ to denote the original left currents
$J_n$.     
 
\begin{prop} {\em (Right-currents on the lattice)}  With
$J^{\R}_n \in \G_a \o \K_N$ defined as above, one finds that
\begin{enumerate}
\item the elements $J^{\R}_n$ and $J^{\L}_m$ commute for arbitrary $n,m$, 
    $$   \Je{1}{\R}_n \Je{2}{\L}_m = \Je{2}{\L}_m \Je{1}{\R}_n\ \ ; $$
\item the elements $J^{\R}_n$ satisfy the following exchange and 
   functoriality relations
   \ba 
    \Je{2}{\R}_{n+1} R \Je{1}{\R}_n & = & 
     \Je{1}{\R}_n \Je{2}{\R}_{n+1} \ \ , \\[1mm]
     \Je{2}{\R}_n \Je{1}{\R}_n
     &  = & R \D_a (J^{\R}_n)  
    \ \ . 
   \ea
(Here we are using the same notations as in the definition of the 
lattice current algebra above.)
\end{enumerate}
\end{prop}

If we denote the subalgebra in $\K_N$ generated by the components 
of left currents $J^{\L}_n$ by ${\cal J}^{\L}_N$ and, similarly, 
use ${\cal J}^{\R}_N$ to denote the subalgebra generated by components
of $J_n^{\R}$, the result of this proposition can be summarized in the 
following statement: ${\cal J}^{\R}_N$ and ${\cal J}^{\L}_N$ form commuting 
subalgebras in $\K_N$, and ${\cal J}^{\R}_N$ is isomorphic to 
$({\cal J}^{\L}_N)_{op}$. Here the subscript $\ _{op}$ means 
opposite multiplication. Both statements have their obvious 
counterparts in the continuum theory (cp. Eq. (\ref{Rcont})).  
\vspace*{5mm}

\noindent
{\sc Example:} {\em (The right $U(1)$-currents)}  
We continue the discussion of the $U(1)$-current algebra 
on  the lattice by  constructing the right currents. Our 
general theory teaches us to consider
$$     w^{\R}_n = g_{n} w_n^{-1} g_{n-1} \ \ . $$
Here $g_n = N^1_{n,+} = \iota_n  (g)$ and  $g_{n-1} = 
N^1_{n-1,-}=  \iota_{n-1}(g)$. The reader is invited to 
check that these elements commute 
with $w^{\L}_n = w_n$. 

\ssection{Monodromies}
\def\cC{{\cal C}}
In the continuum theory one is particularly interested in the 
behavior of the chiral fields $g^\cC (x) = P \exp(\int_0^x j^\cC (x)dx)$ 
under rotations by $2\pi$, i.e. the {\em monodromy of $g^\cC$}. 
Here and in the following, $\cC$ stands for either $\L$ or $\R$. 
The monodromy of $g^\cC$ is determined by the expression
$$ m^\cC = P \exp (\oint j^\cC (x) dx). $$ 
Due to 
the regularizing effect of the lattice, left and right monodromies, 
$M^{\L}, M^{\R}$, are relatively easy to construct and control for  
our discrete current algebra. They are obtained as an ordered product 
of the chiral {\em lattice holonomies} $J^{\L}_n$ or $J^{\R}_n$ along
the whole circle, i.e. 
\ba
    M^{\L} & = &  v_a^{1-N} J^{\L}_1 J^{\L}_2 \cdots J^{\L}_N \ \  
    \ \ \ \mbox{and}
    \\[1mm]
    M^{\R} & = & v_a^{1-N}  J^{\R}_N \cdots J^{\R}_2 J^{\R}_1 \ \ . 
\ea
When we derive relations for the monodromies, it is convenient to 
include the factors involving $v_a = v \o e \in \G_a \o \K_N$.  
The definition in terms of left and right currents produces elements
$M^{\L}$ and $M^{\R}$ in $\G_a \o \K_N$. Their algebraic structure 
differs drastically from the properties of the currents $J_n$. We 
encourage the reader to verify the following list of equations:
\ba 
\M{1}{\L} R \M{2}{\L} =  R \D_a (M^{\L}) \ \ & , & \ \ 
\M{2}{\R} R \M{1}{\R} =  R \D_a (M^{\R}) \ \ , \nn \\[1mm] 
\D_0(\xi) M^{\L} = M^{\L} \D_0(\xi) \ \ & , & \  \ 
\D'_0(\xi) M^{\R} = M^{\R} \D'_0 (\xi) \  \  \nn
\ea
for all  $\xi\in \G$ and with $\D'_0(\xi) = (id \o \iota_0)(\D'(\xi))
=  (id \o \iota_0)(\s  \circ \D (\xi))$. The functoriality relations
are familiar from Subsection 2.2 and imply that the algebra generated 
by components of the monodromy $M^{\L}$ or $M^{\R}$ is isomorphic 
to $\G$ or $\G_{op}$. There are several places throughout the 
paper where this observation becomes relevant for a better 
understanding of our results. 

As usual, we may  evaluate the elements $M^\cC$ in 
irreducible representations of $\G_a$. This results in a set 
of $\K_N$ valued matrices $M^{\cC,I} \equiv  (\t^I \o id) (M^\cC)$. 
Their quantum traces  
$$        c^I_\cC \equiv \tr^I_q(M^{\cC,I}) \ \  $$
are  elements in the algebra $\K_N$ which have a number of  
remarkable properties. First, they are central elements in  
the lattice current algebra $\K_N$, i.e. the $c^I_\cC$  
commute with all elements $A \in \K_N$. Even more important
is that $c^I_\R, c^I_\L \in \K_N$ generate two commuting       
copies of the {\em Verlinde algebra} \cite{Ver}. Explicitly 
this means that  
$$     
      c^I_\cC c^J_\cC = \sum N^{IJ}_K c^K_\cC \ \ \mbox{ and } 
       \ \ (c^K_\cC)^* = c^{\bar K}_\cC \ \  
$$
for $\cC = \R,\L$. Here and in the following $\bar K$ denotes the 
unique label such that $N^{K \bar K}_0 = 1, $ and $0$ stands for 
the trivial representation $\t^0 = \e$. If the $S$-matrix  
$S_{IJ} = \N (\tr^I_q \o \tr^J_q)(R'R)$ is invertible, and 
$\N$ is suitably chosen, the linear combinations 
$$   \chi_\cC^I = \sum_J \N d_I S_{I \bar J}c^J_\cC $$
provide a set of orthogonal central projectors in $\K_N$, for each 
chirality $\cC = \R,\L$, i.e.,
$$   \chi^I_\cC \chi^J_\cC = \d_{IJ} \chi^I_\cC \ \ , \ \ 
     (\chi^I_\cC )^* = \chi^I_\cC $$
Proofs of all these statements can be found in 
\cite{AGS}. We will see in the next section that products 
$\chi^I_\L \chi^J_\R$ provide a complete set of minimal central 
projectors in the lattice current algebra or, in other words, 
they furnish a complete set of characteristic projectors 
for the irreducible representations of $\K_N$. 
\vspace*{5mm}

\noindent
{\sc Example:} {\em  (The center of the lattice $U(1)$ 
current algebra)} In terms of the variable $w_n =  q^{(1-p)/2}
J^1_n = -q^{1/2} J^1_n$ (cf. Subsection 2.3 for notations), 
the definition  of the monodromy $M^{\L,1} = (\t^1 \o id)  
(M^\L)$ becomes 
$$   M^{\L,1} = (-1)^N q^{3N/2-1} w_1 w_2 \cdots w_N 
     \ \ \in \K_N\ \ .   $$
In this particular example, the quantum trace is trivial  
so  that $c^1_\L = M^{\L,1}$. It is easily checked that 
$c^1_\L$ commutes with all the generators $w_n, g_n \in  
\K_N$ and that it satisfies 
$$       (c^1_\L)^p = c^0_\L = e \ \ . $$
Of course, this relation follows also from the 
formula $c^r_\L c^s_\L = \sum N^{rs}_t c^t_\L = c^{r+s}_\L$.     
Similar considerations apply to the right currents.  

\ssection{The inductive limit $\K_\infty$} 
So far, the lattice current algebras $\K_\N$ depend on the 
number $N$ of lattice sites, and one may ask what happens when 
$N$ tends to infinity. A mathematically precise meaning 
to this question is provided by the notion of {\it inductive limit}.
The latter requires an explicit choice of embeddings of 
lattice current algebras for different numbers of lattice 
sites. They will come from some kind of inverse block-spin
transformation \cite{FaGa}. 

Suppose we are given two lattice current algebras $\K_N$ and 
$\K_{N+1}$ with generators $J_n, n=1,\dots, N$ and $\hat J_m, m = 
1,\dots, N+1$ respectively. The embeddings of $\G$ into $\K_N$
or $\K_{N+1}$ will be denoted by $\iota_n$ or $\hat\iota_m$. An 
embedding $\c_N: \K_N \rightarrow \K_{N+1}$ is furnished by 
\ba
   \c_N (J_n) & = & \hat J_n    \mbox{ for all }  \ \ n < N\nn,  \\[1mm]
   \c_N (J_N) & = & v_a^{-1}\hat J_{N} \hat J_{N+1} \ \ 
   \mbox{ and } \nn \\[1mm]
   \c_N (\iota_n(\xi)) & = & \hat \iota_n(\xi) \ \ \mbox{ for all } 
       n < N \ \ , \nn \\[1mm]
   \c_N(\iota_N(\xi)) & = & \hat \iota_{N+1} (\xi)\ \ . \nn 
\ea    
The intuitive idea behind $\c_N$ is to pass from $\K_N$ to 
$\K_{N+1}$ by dividing the $N^{th}$ link on the lattice of 
length $N$ into two new links, so that we end up with a lattice
of length $N+1$. Observe that $\c_N$ maps the monodromies 
$M \in \G_a \o \K_N$ to $\hat M \in \G_a \o \K_{N+1}$, and, consequently, 
the same holds  for our projectors $\chi^I_\L, \chi^J_\R \in \K_N$
and $\hat \chi^I_\L, \hat \chi^J_\R \in \K_{N+1}$,   
\be 
     \c_N (\chi^I_\L) = \hat \chi^I_\L \ \ , \ \ 
    \c_N (\chi^J_\R) = \hat \chi^J_\R\ \ . 
\label{chiemb} 
\ee
Since the set of numbers $N$ is directed, the collection of $\K_N$,
together with the maps $\c_N$, forms a directed system, and we can 
define the inductive limit 
$$ \K_\infty \equiv  \lim_{N \rightarrow \infty} \K_N \ \ . $$
By definition, $\K_\infty = \bigcup_N \K_N / \sim$ where two elements 
$A_N \in \K_N$ and $A_{N'} \in \K_{N'}$ are equivalent, iff
$A_N$ is mapped to $A_{N'}$ by a string of embeddings $\c_M$, 
i.e., $A_{N'} = \c_{N'-1} \circ \dots \circ \c_{N+1}\circ \c_N(A_N)$. 
For the lattice $U(1)$-current algebra, a detailed 
investigation of this inductive limit was performed in 
\cite{AlRe}. 
 
We have chosen to define the block-spin transformation by dividing
the $N^{th}$ link of the lattice. Now we introduce another block
spin operation by dividing the $1^{st}$ link of the lattice:
\ba
   \cc_N (J_n) & = & \hat J_{n+1}    \mbox{ for all }  \ \ n > 1\nn,  \\[1mm]
   \cc_N (J_1) & = & v_a^{-1}\hat J_{1} \hat J_{2} \ \ 
   \mbox{ and } \nn \\[1mm]
   \cc_N (\iota_n(\xi)) & = & \hat \iota_{n+1}(\xi) \ \ \mbox{ for all } 
       n  \ \ . \nn 
\ea    
Notice that the two block-spin operations `commute' with each other:
$\cc_{N+1} \circ \c_{N} = \c_{N+1} \circ \cc_N$. While the map $\c_N$
is used in the definition of the inductive limit, we reserve $\cc_N$
for the definition of the co-product for lattice current algebras
(see Section 4).

\section{Representations of the Lattice Current Algebra} 
\setcounter{equation}{0}
\def\cU{{\cal U}}
 
The stage is now set to describe our main result on the 
representation theory of the lattice current algebras. 
We will begin with a much simpler problem of representing 
two important subalgebras of $\K_N$. Their representations
will serve as building blocks for the representation theory
of the lattice current algebra $\K_N$.     

\ssection{The algebra ${\cal U}$} 
The lattice current algebras $\K_N$ contain $N-1$ non-commuting 
{\em holonomies} $U_\nu, \nu = 1, \dots, N-1,$ 
$$   U_\nu \equiv  v_a^{1-\nu} J_1 \cdots J_\nu \ 
     \in \G_a \o \K_N\ \ .  $$
The objects $U_\nu$ commute with all elements in the symmetry
algebras  $\G_n$, except the ones for $n = 0$ and $ n = \nu$. These
properties of $U_\nu$ remind us of holonomies in a gauge 
theory, which transform nontrivially only under gauge 
transformations acting at the endpoints of the paths. Together, the 
elements  in $\G_0 \o \G_\nu \subset  \K_N$ and the components 
of $U_\nu$  generate a subalgebra, ${\cal U}_\nu$, of the 
lattice current algebra $\K_N$. These  subalgebras $\cU_\nu 
\subset \K_N$ are all isomorphic to the algebra $\cU$ which 
we investigate in this subsection. We begin with a  
definition of $\cU$.     

\begin{defn}  The algebra 
${\cal U}$ is the *-algebra generated by components of 
elements $U, U^{-1} 
\in \G_a \o {\cal U}$ together with elements in $\G_0 \o 
\G_1$ such that 
\ba
    \U{1}{} \U{2}{}  & = & R \D_a(U)
      \ \ , \nn\\[1mm]
   \iota_1(\xi) U   =   U \D_1(\xi) \ \ & , & 
    \D_0(\xi) U = U \iota_0 (\xi) \nn
\ea
and $U^{-1}$ is the inverse of $U$. Here we use the same notations
as in Subsection 2.3. In particular, $\iota_{0,1}$ denote the 
canonical  embeddings of $\G$ into $\G_0 \o \G_1$. The $*$-operation
on ${\cal U}$ extends the $*$-operation on $\G_0 \o \G_1
\subset \cU$, so that
 $$     U^* = S_1^{-1} U^{-1}  S_0\  \ . $$
Here $S_i  = (id \o  \iota_i) (\D(\k^{-1}) (\k  \o \k)R^{-1}) \in 
\G_a \o  \cU$ for $i = 0,1$, and $\k$ is a certain central square root of 
the ribbon element, as  before.    
\end{defn}
 
The algebra $\cU$ admits a very nice irreducible representation, 
$D$, which is constructed by acting with elements in ${\cal U}$ 
on a ``ground state'' $\vac$. The state $\vac$ may be characterized 
by the following (invariance-) properties
$$    
    \iota_i(\xi)  \vac = \vac \e(\xi)    \ \
     \mbox{ for all } \ \  \xi  \in \G \ \  \ \ i = 0,1 
$$
where $\e$ is the trivial representation (co-unit) of $\G$. 
Here and in the following we neglect to write the letter
$D$ when elements in $\cU$ act on vectors. Since we are dealing 
with  a unique representation of $\cU$, ambiguities are  
excluded. While the preceding formula means that $\vac$ is 
invariant under the action of 
$\G_0 \o \G_1$, the components of $U \in \G_a \o \cU$ create new 
states in the carrier space, $\Re$, of the representation $D$, 
$$  |\pi \rangle \equiv \pi (U) \vac = (\pi \o id )(U) \vac \ 
      \ \in \Re $$
for all  $\pi \in \G'_a$. In particular, one identifies $| \e\rangle
= \vac$ because $\e(U) =  U^0 = e$. We can think of the vectors 
$|\pi \rangle$ as coming from a universal object $u \in  \G_a \o \Re$, 
i.e., 
$$   |\pi \rangle \equiv \pi(u) = (\pi \o id)(u) \ \ . $$
In other  words,  $u =  U \vac$. 
A complete description  of the representation $D$ on $\Re$ 
is given in the following proposition. 

\begin{prop} {\em (Representation of ${\cal U}$)} There exists  
an irreducible *-represen\-tation $D$ of the algebra ${\cal U}$ 
on a carrier space $\Re$ such that \label{Urep}
\ba    
   \U{1}{} \uE{2}{} & = & R \D_a(u) 
   \ \ , \nn \\[1mm]
  (e\o  \iota_1(\xi)) u  =  u (\xi \o e)
 \ \ \ & , & \ \ (e \o \iota_0(\xi)) u = (\S(\xi) \o e) u  
  \ \ . \nn
\ea
Here $u \in \G_a \o \Re$ and $\D_a (u) = (\D \o id)(u)$. The 
representation space $\Re$ contains a unique invariant 
vector $\vac \in \Re$. 
\end{prop}

{\sc Proof:} The formulas for the action of ${\cal U}$ on 
$\Re$ follow from $u = U \vac \in \G_a \o \Re$ by using the  
invariance of $\vac$ under the action of $\iota_1(\xi)$ and 
$\iota_0(\xi)$. In particular, we have that 
\ba 
    \U{1}{} \uE{2}{}  &  = & 
    \U{1}{} \U{2}{}\vac = R \D_a(U) \vac \nn \\[1mm]
     & = &  R \D_a(u) \ \  \  \mbox{ and } \nn  \\[2mm]
    (e\o  \iota_1(\xi)) u & = &  U  \D_1 (\xi) \vac = 
       U \vac (id \o \e)(\D(\xi))  \nn \\[1mm]
      & = & u (\xi \o e)\ \ . \nn 
\ea 
To derive the action of $\G_0$ on $\Re$ we rewrite the equation
$U \iota_0(\xi) = \D_0(\xi) U$ according to  
$$   
   \iota_0(\xi) U = ( \S(\xi^1_\vs)\o e) U 
     (e \o \iota_0(\xi^2_\vs))
      \ \ . $$
Here we have inserted the expansion $\D(\xi) = \xi^1_\vs \o \xi^2_\vs$ 
and used some standard Hopf-algebra properties of the antipode $\S$.     
Then one proceeds as in the computation of $\iota_1(\xi) u$
to obtain the last formula claimed in Proposition \ref{Urep}.

Observe that the formulas in the preceding 
Proposition  define an action -- not a co-action -- of the 
algebra ${\cal  U}$ on the representations space $\Re$. 
Again, it is important to keep in  mind that the co-product 
$\D$ in the first formula of Proposition \ref{Urep} acts on the 
first  tensor factor $\G_a$  of  $u \in \G_a \o \Re$. In terms 
of the multiplication $\ast$ in $\G'$ (cf. Subsection 2.3)
one has that 
$$   \pi_1(U) |\pi_2\rangle = |\pi_1 \ast \pi_2 \rangle $$
for all $\pi_1, \pi_2 \in \G'$. Consequently, the components 
$\pi(U)$ act on $\Re$ as some kind  of (twisted) multiplication 
operators. Furthermore, the matrix elements of $u^I =  \t^I (u)    
\in \End(V^I) \o \Re$ span a $\d_I^2$-dimensional subspace of 
$\Re$ which is invariant under the action of $\G_0 \o \G_1  
\subset  \cU$. With the help  of Proposition \ref{Urep} we 
deduce
$$  \iota_1(\xi) u^I =  u^I \t^I(\xi) \ \ , \ \ 
    \iota_0(\xi) u^I =  \t^I(\S(\xi)) u^I $$
for all $\xi \in \G$ and with $\iota_i(\xi) = 
(e^I \o \iota_i(\xi))$, where $e^I$ is the unit element in 
$\End(V^I)$. The two formulas furnish the following decomposition 
of $\Re$ into a direct sum of $\G_0 \o \G_1$-modules,    
\be   
     \Re \cong 
     \bigoplus_I \check{V}^I \o V^I \  \ , \label{dec} 
\ee
where $\check{V}^I$ is dual to $V^I$ and the sum extends over 
the classes of irreducible representations of $\G$. $\check{V}^0 
\o V^0  \subset \Re$ coincides with the one-dimensional subspace 
spanned by $\vac$. All these features of the representation 
space $\Re$ resemble those of the algebra of square-integrable 
functions on a Lie-group
$G$ with its characteristic action of left and right invariant
vector fields. This  similarity is not too surprising and 
can be traced back to the analogy between the objects $U$ and 
$T$. In Subsection 2.3, $T$ was found to satisfy  ``$RTT$-
relations'' which are the key ingredient in the deformation 
theory of groups $G$.

Mainly for technical reasons we finally 
look at a certain subalgebra ${\cal D}$ of $\cU$ and its  
action on $\Re$. 
 
\begin{lemma} Let ${\cal D}$ be the subalgebra of $\cU$ which is 
generated by elements $\xi \in \G_1$  and components of $U$.  
When the representation $D$ of ${\cal U}$ is restricted to 
${\cal D} \subset \cU$ it furnishes an irreducible representation
$D$ of ${\cal  D}$ on $\Re$. \label{Dirred}
\end{lemma} 

{\sc Proof:} 
To prove this Lemma we show that every vector $|\pi\rangle$  in the 
representation space $\Re$ is cyclic under the action of ${\cal D}$. 
Since $\Re$ contains the cyclic vector $\vac$, our task simplifies 
to the following problem: show that for every $|\pi\rangle \in \Re$ 
there is a representation operator $A_\pi\in D({\cal D})$ such 
that $A_\pi  |\pi\rangle = \vac$. For the proof it is crucial to 
find the projector on $\vac$ in $D({\cal D})$. It is constructed   
from the minimal central projector $P^0 \in \G$ that corresponds 
to the trivial representation $\e = \t^0$. By definition, 
$P^0$ satisfies $\t^I(P^0) = \d_{I,0}$, so that 
$$    \iota_1( P^0)  | \pi' \rangle   = \vac \pi'(P^0)   \ \  $$
holds for all $\pi' \in \G'$. The other ingredient we  need below 
is a distinguished element $\mu \in \G'$ -- called the {\em  right 
integral} of $\G$  --  with the properties  
$$   (\mu \o id) \D(\xi) =  \mu(\xi) \ \ ,\ \ \mu (P^0) =  1 \ \ . $$
Now we choose $\xi_\pi \in \G$ such that $\pi (\xi_\pi)
= 1 $. A short  technical computation shows that
$$   (e \o \xi_0) \D(P^0) = (\S(u^{-1} \xi_0 ) \o  e) R \D(P^0) \ \ ,  $$
where $u = g^{-1} v \in \G$ and  $g$ was introduced in Subsection 2.1. 
We abbreviate 
$\eta_0 \equiv \S(u^{-1} \xi_0 )$ and regard $\eta_0$ as a map from  
$\G$ to $\G$ acting by left multiplication so that $\mu \circ \eta_0$
makes sense as an  element in $\G'$. Let us  define 
$$        A_\pi \equiv \iota_1(P^0) (\mu \circ \eta_0)(U) \ \ .  $$
The following calculation proves that $A_\pi |\pi \rangle = \vac$ 
and hence completes the proof of the lemma.           
\ba
       A_\pi |\pi  \rangle  
       & = & \iota_1( P_0) (\mu\circ \eta_0)(U) |\pi\rangle 
                             \nn \\[1mm]
       & = & \iota_1( P_0) |(\mu\circ \eta_0) \ast \pi\rangle 
                             \nn \\[1mm]
       & = & \vac ((\mu\circ \eta_0) \ast \pi)(P_0)
                             \nn \\[1mm]
       & = & \vac (\mu\o \pi)((\eta_0 \o e) R\D(  P_0))
                             \nn \\[1mm]
       & = & \vac  \pi(\xi_0) \mu(P^0) = \vac\ \ . 
                             \nn 
\ea
The algebra ${\cal  D}$ is a semidirect product of $\G$ and its dual $\G'$, 
the latter being  supplied with the twisted product $\ast$ that we 
discussed in Subsection 2.3. In this light, ${\cal D}$ appears as a close
relative of the deformed  cotangent bundle $T_q^*G$ over a group $G$ 
which differs from the structure of ${\cal D}$ only through the use of the 
standard product $\cdot$ in $\G'$. 
\vspace*{5mm}

\noindent
{\sc Example:} Let us  continue our tradition and  illustrate the 
theory with  the example $\G = \Z_q$. The algebra $\cU$ 
is then generated by the unitary elements $g_0,g_1$ and $w$ satisfying
  $$   w^p = 1 \ \ , \  \   g_1 w = q w g_1 \ \ , 
      \ \  w g_0 =qg_0 w $$
and $g_0$ commutes with $g_1$ (cf. Subsection 2.3 for further 
details). States in $\Re$ are created from a ground state 
$\vac$  with invariance properties 
$$    g_1 \vac = \vac \ \ \mbox{ and }\ \ g_0 \vac = \vac\\ .  $$
Through iterated application of $w$ on $\vac$ we may produce 
$p$ linearly independent vectors 
$$     | r \rangle \equiv  w^r \vac \  \in \Re \ \ \mbox{  for } 
   \ \ r = 0, \dots, p-1 \ \ . $$
Specializing the proof of Proposition \ref{Urep} to our example, we obtain 
$$    g_1 |r \rangle  = q^r w^r g_1 \vac = q^r |r\rangle  $$
and similarly for $g_0$. This shows that the one-dimensional 
subspace spanned by $|r \rangle$ corresponds to  the 
summand $\check{V}^r \o V^r$ in the decomposition (\ref{dec}) 
of  $\Re$. The subalgebra ${\cal D}$ of $\cU$ is generated by 
$g = g_1$ and $w$,  with Weyl  commutation relations 
$$    w g  = q   g w   \ \  $$
and $g^p = e =  w^p$. It acts irreducibly on the $p$-dimensional 
space $\Re$.

\ssection{The algebra $\K$}
Before we deal with the general situation, it is helpful 
to study the simplest example of a lattice current algebra for 
which the lattice consists of only one (closed) link and one site, 
i.e., the case $N=1$. Strictly speaking, $\K_1$ has not been defined 
above. So we must first give a definition.         

\begin{defn} {\em (The algebra $\K$)} The *-algebra $\K \equiv \K_1$ 
is generated by components of $M, M^{-1} \in \G_a \o \K$ and elements 
$\xi \in \G$ with the following relations
\ba 
     \M{1}{} R \M{2}{} & =& R \D_a(M)
     \ \ ,\nn \\[1mm]
    \D(\xi) M & = & M \D(\xi) \ \ \mbox{ for all } \ \ 
     \xi \in \G \ \  \nn
\ea
$M^{-1}$  is the inverse of $M$, so that $M^{-1} M = e = M M^{-1}$. 
The action of $*$ is extended from $\G$ to $\K$ by the formula
$$   M^*  = S^{-1} M^{-1} S\ \ , $$
where  $S = \D(\k^{-1}) (\k \o \k) R^{-1}$,  as before. 
\end{defn}
 
Components $\pi (M)$ of the monodromy $M$ can be represented 
on the carrier spaces $V^I$ of the representations $\t^I$ of 
$\G$. This is accomplished by the formula
$$  D^I( \pi(M)) = (\pi \o \t^I)(R'R) \ \ \in \End (V^I)$$ 
for all linear forms $\pi \in  \G'_a$ on $\G_a$. An equivalent  
universal formulation without reference to linear 
forms $\pi$ is
$$  D^I(M) =  (id \o \t^I)(R'R)\  \in \G_a \o \End(V^I)\ \ . $$ 
Indeed, one may check that such an action on $V^I$ is consistent 
with the functoriality of $M$, i.e., with the first relation in 
Definition 5,
\ba  D^I(\M{1}{} R \M{2}{}) & = & 
      (id \o  id \o \t^I)[ R'_{13} R_{13} R_{12} 
        R'_{23} R_{23} ] \nn \\[1mm]
   &=& (id \o  id \o \t^I)[ R'_{13} R'_{23} 
        R_{12} R_{13} R_{23} ] \nn \\[1mm]
   &=& (id \o  id \o \t^I)[ R_{12} R'_{23} 
        R'_{13} R_{13} R_{23} ] \nn \\[1mm]
   &=& (id \o  id \o \t^I)[ R_{12} 
         (\D \o id) (R' R)] \nn \\[1mm]
   &=& R (\D_a \o id) (D^I(M)) \nn \ \ . 
\ea
{}From our discussion in Subsections 2.4, 2.5 we know already that 
$\G_a \o \K$ contains not only the left monodromy $M^\L =  M$ but also 
the right monodromy 
\ba    M^\R & = &  N^{-1}_- M^\L N_+  \ \in \G_a \o \K \ \  
       \  \mbox{ with }  \nn \\[1mm]   
     N_+  = R'  \ \  \ & , & \  \  N_-  = R^{-1} \ \ . \nn 
\ea
Here $N_\pm$ are regarded as elements in $\G_a \o \G \subset \G_a \o \K$,
as before. Since the components of $M^\L$ and  $M^\R$ commute, they  can  be
represented on spaces $V^I \o V^J$ such that the right/left-monodromies 
act trivially on the first/second  tensor factor, respectively. This action 
of left and right monodromies on 
$$ W^{IJ} \equiv V^I \o V^J \ \  $$
can be extended to an action of the entire algebra $\K$. 
Actually, the algebra $\K$ is isomorphic to the Drinfeld
double of $\G$ \cite{ReSTS,FN}.

\begin{prop} {\em (Representations of $\K$)} \label{Krep}
The algebra $\K$ has a series of irreducible 
*-representations, $D^{IJ}$, defined on the spaces $W^{IJ}$. 
Explicitly, the action is given 
by 
\ba 
     D^{IJ} ( M) & = & (id \o \t^I) (R'R) \ \  , \nn\\[1mm]
  D^{IJ}(\xi) & = &  (\t^I  \bo \t^J)(\xi) = 
       (\t^I \o \t^J) (\D(\xi)) \ \ .\nn
\ea
Here $D^{IJ} (M) \in \G_a \o  \End (V^I)$ is regarded as an 
element of $\G_a \o \End(W^{IJ})$ with trivial action  on
the second tensor factor $V^J$ in  $W^{IJ}$. 
\end{prop}     

{\sc Proof:} To check the representation property is left as an 
exercise to the reader. It may be helpful to consult
Theorem 12 of \cite{AlSc}.
Irreducibility  follows from the fact that the action  of the 
monodromies on the  spaces $V^I$ is irreducible (cf. Lemma 1 of 
\cite{AlSc}).  
\vspace*{5mm}

\noindent
{\sc Example:} For  $\G = \Z_q$, the definition 
of $\K$ furnishes an algebra with generators $c = c^1_\L$ and 
$g$ such that 
$$   c^p = 1 \ \ , \ \  g^p = 1  \ \ \mbox{ and } 
     c g  =  g c \ \ . $$
Some explanation can  be found at the end of Subsection 2.5.
The abelian algebra $\K$ has one-dimensional representations,
$D^{st}$, on spaces $W^{st}$ labeled by two integers $s,t = 0, \dots, 
p-1$,  
$$  D^{st}(c) = q^s \ \ \ , \ \ \  D^{st}(g) = q^{s+t} \ \ ,$$ 
where $q = exp(2\pi i/p)$. For comparison with the general 
formulas in Proposition \ref{Krep}, one should keep in mind that 
$g$ is a generator of the algebra $\G$ while $c$ coincides
with the element $\t^1(M) \in \K$, up to a scalar factor.            

\ssection{Representations of $\K_N$}
In representing the full lattice current algebras $\K_N$ it is 
convenient to pass to a new set of generators: Let us define
elements $U_\nu \in \G_a \o \K_N, \nu = 1, \dots, N-1,$ as above by 
$$ U_\nu \equiv  v_a^{1-\nu} J_1 \cdots J_{\nu}\ \ . $$
Then the $U_\nu$, together with $M \equiv M^\L$ and the elements
$\xi \in \G_n, n = 1, \dots ,N,$ generate $\K_N$.  They obey 
the following relations 
\ba
    \U{1}{}_\nu \U{2}{}_\nu   & = & R \D_a( U_\nu) 
    \ \ , \nn \\[1mm]
    \M{1}{} R \M{2}{} & = & R  \D_a( M)
     \ \ , \nn \\[1mm]
    R' \U{1}{}_\nu  \U{2}{}_\mu  &=& 
     \U{2}{}_\mu \U{1}{}_\nu \ \ \mbox{ for }  \ 1 \leq \nu < \mu \leq
      N-1 \ \ , \nn \\[1mm]
    R' \U{1}{}_\nu \M{2}{} &  = & \M{2}{} 
     R'\U{1}{}_\nu \ \ , \\[2mm]
    \D_0(\xi) M & = & M \D_0 (\xi) \ \ \mbox{ for }\ \  
      \xi  \in \G \ \ ,\nn \\[2mm]
    \iota_\nu(\xi) U_\nu & = & U_\nu \D_\nu (\xi) \ \ \mbox{ for }\ \  
      \xi  \in \G  \ \ , \nn\\[2mm]
    \D_0(\xi) U_\nu & = & U_\nu \iota_0(\xi) \ \ \mbox{ for }\ \  
      \xi  \in \G \ \ . \nn
\ea
We infer that $\K_N$ contains $(N-1)$ copies of the 
algebra ${\cal U}$ generated by $U_\nu, \G_\nu, \G_0$ and one copy of 
the algebra $\K$ generated by $M \equiv M^\L$ and $\G_0$. 
These subalgebras do not commute, but the non-commutativity is felt 
only by $U_\nu, M $  and $\G_0$. 
In any case, the exchange relations motivate us to look for 
representations of the lattice current algebra $\K_N$ on 
spaces 
\be    W^{IJ}_N \equiv \Re^{\o_{N-1}} \o W^{IJ} \ \ .\label{repsp} \ee 
To state the formulas we introduce the notation $D_\nu, \nu = 1, 
\dots, N-1,$ which stands for the representation 
$D$ of the algebra ${\cal U}_\nu$ 
on the $\nu^{th}$ factor of the tensor product (\ref{repsp}), i.e., 
for every $X \in {\cal U}_\nu $, 
$$  D_\nu(X) = id_1 \o \dots \o id_{\nu-1} \o D (X) \o id_{\nu+1} \o 
\dots \o id_{N-1} \o Id\ \ , $$
where $id_\mu$ acts as the identity on the $\mu^{th}$ factor $\Re$ in 
$W^{IJ}_N$, and $Id$ is the identity on $W^{IJ}$. Similarly, $D^{IJ}$
denotes the action of the subalgebra $\K \subset \K_N$ on the 
last factor $W^{IJ}$ in $W^{IJ}_N$. 

The map $\iota_0$ embeds elements of $\G$ into all the subalgebras 
$\cU_\nu$ of the lattice current algebra $\K_N$. Hence $\G$
acts on each tensor factor $\Re$ in $W^{IJ}_N$ independently with the 
help of representations $D_\nu$. We employ the co-product 
$\D$ on $\G$ to obtain a family $\vth_\nu, \nu = 1, \dots, N,$ of 
$\G$-actions on $W^{IJ}_N$. In the following, $\vth_1$ is the trivial 
representation, $\vth_1 = \e$, and 
$$ 
    \vth_\nu (\xi) =  (D_1 \bo \dots  \bo  D_{\nu-1})(\iota_0(\xi)) 
                   \o id_{\nu} \o \dots id_{N-1}\o Id  
$$
for all $\xi \in  \G$. The symbol $\bo$ denotes the tensor product 
of representations. 

With these conventions we are prepared to define 
representations of $\K_N$ on $W^{IJ}_N$. The essential idea is 
borrowed from the well known Jordan-Wigner transformation. 
In fact, in writing the actions of our generators on $W^{IJ}_N$  
we have to relate the  $U_\nu$'s and $M$ (which obey non-trivial
exchange relations among each other) to the operators $D_\nu
(U_\nu)$ and $D^{IJ}(M)$. The latter act on different tensor 
factors in $W ^{IJ}_N$ and hence commute. In  analogy  with  
the `tail'-factors  $\prod_{i=0}^n \s^3_i$ of the Jordan-Wigner
transformation, we will employ tails of  $R$-matrices to  
express the $U_\nu$ in terms of $D_\nu(U_\nu)$ and $M$ in 
terms of $D^{IJ}(M)$. Related constructions appear  in Majid's 
`transmutation theory' (see e.g. \cite{Maj}).   
 
\begin{theo} {\em (Representations of $\K_N$)} The algebra $\K_N$ 
has a series of irreducible *-representations $D_N^{IJ}$ realized
on the spaces $W^{IJ}_N$ given in equation (\ref{repsp}). 
In the representation $D^{IJ}_N$, the generators of $\K_N$ act as
\ba 
    D^{IJ}_N (U_\nu) & = & (id \o \vth_\nu) (R^{-1}) 
                D_\nu( U_\nu)\ \ , \nn \\[1mm]
    D^{IJ}_N (M) & = & (id \o \vth_{N}) (R^{-1})
                D^{IJ} (M) (id \o \vth_{N}) (R) \ \ , \nn\\[1mm]
    D^{IJ}_N (\iota_\nu(\xi)) & = &  D_\nu (\iota_\nu(\xi))\ \ \mbox
        {for } \ \ \nu=1,\dots,N-1\  ,\ \ \xi \in \G \ \ , \nn\\[1mm]
    D^{IJ}_N (\iota_0(\xi)) & = & (\vth_{N} \o D^{IJ}) (\D_0(\xi)) 
     \ \ \mbox{ for all } \ \ \xi \in \G\ \ .  \nn
\ea
Every *-representation of $\K_N$ can be decomposed into a 
direct sum of the irreducible representations $D^{IJ}_N$. 
\end{theo}
 
{\sc Proof:} The proof is similar to the one of Theorem 15 in 
\cite{AlSc}. Although we do not present the details, 
we stress that the factors $(id \o \vth_\nu)
(R^{\pm 1})$ produce ``tails of $R$-elements'' which are responsible 
for the correct exchange relations of $U_\nu, M$ in the representations 
$D^{IJ}_N$. From the definition of $\vth_\nu$ and quasitriangularity of 
$R$ we infer (for $\nu > 1$)
$$   (id \o \vth_\nu)(R^{-1}) = R^{-1}_{12} R^{-1}_{13} \dots
        R^{-1}_{1\nu} \  \ . $$
Here it is understood that the second tensor factor of $R_{1\mu}$ 
is represented on the $\mu^{th}$ factor $\Re$ in $W^{IJ}_N$ with  
the help of  $D_\nu \circ \iota_0$. A simple  example suffices 
to illustrate how nontrivial exchange relations between $U_\nu,M$ 
arise in our representations
\ba 
   D^{IJ}_N( R'_{12}  \U{1}{}_1 \U{2}{}_2) & = & 
    R'_{12} D_1(\U{1}{}_1) R^{-1}_{23} D_2(\U{2}{}_2) \nn \\[1mm]
 & = &    R'_{12} (R'_{12})^{-1} R_{23}^{-1} 
    D_1(\U{1}{}_1) D_2(\U{2}{}_2) \nn \\[1mm]
 & = &    R_{23}^{-1} D_2(\U{2}{}_2) D_1(\U{1}{}_1) \nn \\[1mm] 
  & = &   D^{IJ}_N (\U{2}{}_2  \U{1}{}_1 ) \ \ . \nn
\ea  
These equations are to be understood as equations in 
$\G_a \o \G_a \o \End(W^{IJ}_N)$. To 
reach  the second line,  we insert the exchange  relation of 
$U_1$ with elements $\iota_0(\xi)$ and quasi-triangularity of 
$R^{-1}$. Then we use that the images of $D_1$ and $D_2$ 
commute. 
 
Irreducibility of the representations $D^{IJ}_N$  follows from 
the case $N=1$ which we have treated in Subsection 3.2, together with 
Lemma \ref{Dirred} of Section 3.1.      
 
It is quite instructive to evaluate the projectors $\chi^K_\L 
\chi^L_\R$ (defined at the end of Subsection 2.5) in the representations
$D^{IJ}_N$. The answer is given by (cf. \cite {AlSc})
$$  D^{IJ}_N( \chi^K_\L \chi^L_\R) = \d_{I,K} \d_{J,L}\ \ . $$
As we have promised above, the elements $\chi^K_\L \chi^L_\R$ are 
characteristic projectors for the irreducible representations
of the lattice current algebra $\K_N$.      

The representation theory described here survives the limit in which
$N$ tends to infinity. With the help of our embeddings $\c_N: 
\K_N \rightarrow \K_{N+1}$ (see Section 2.6) we can define an action 
$D^{IJ}_{N+1} \circ \c_N$ of the lattice current algebra $\K_N$ on 
$W^{IJ}_{N+1}$. It is, of course, no longer irreducible, so 
that $W^{IJ}_{N+1}$ decomposes into a direct sum of irreducible
representations of the algebra $\K_N$. The observations made 
in the preceding paragraph combined with formula (\ref{chiemb})
furnish an isomorphism
$$   W^{IJ}_{N+1}    \cong   W^{IJ}_N \o \Re $$
of $\K_N$ modules, with $\Re$ being a multiplicity space of the 
reduction. In particular, all the irreducible subrepresentations
of the $\K_N$-action on $W^{IJ}_{N+1}$ are isomorphic to $D^{IJ}_N$. 
This implies that the inductive limit for the directed system 
$(\K_N, \c_N)$ splits into independent contributions coming 
from the simple summands of $\K_N$. Since inductive limits of 
simple algebras are simple, we conclude that $\K_\infty$ possesses 
irreducible representations $D^{IJ}_\infty$ on $W^{IJ}_\infty$. 
As for the algebras $\K_N$, the labels $I,J$ run through the 
classes of irreducible representations of $\G$.   
\vspace*{5mm}

\noindent 
{\sc Example:} {\em (Representations of the $U(1)$-current algebra)}
To discuss the representation theory of the $U(1)$-current algebra,
we introduce the generators 
$$ v_\nu = q^{(\nu-1)/2} w_1 \cdots w_\nu \in \K_N $$
which agree with $U^1_\nu = (\t^1\o id)(U_\nu)$, up  to  a scalar 
factor. The whole algebra $\K_N$ is generated from the elements 
$g_\nu \in \G_\nu \subset \K_N$, the  monodromy $c \in \K_N$ and 
the  holonomies $v_\nu \in \K_N$ such that 
\ba    
     g_\nu v_\nu = q v_\nu g_\nu \  \ \  & , & \ \ \ 
     v_\nu g_0 = q g_0 v_\nu \ \ , \nn \\[1mm]
          v_\nu v_\mu = q^{-1} v_\mu v_\nu  & \ & \mbox{ for }  
     \ \nu < \mu \ \  , \nn \\[1mm] 
     g_0^p = g_\nu^p = c^p = v_\nu^p = 1 \ & \ & \mbox{for all  }  \ 
     \nu =  1, \dots, N-1 \nn
\ea
and $c$ commutes with every other element. Vectors  in the 
carrier space of $D^{st}_N$ are denoted by  
$$      |r_1, r_2, \dots,  r_{N-1}\rangle _{s,t}\  \ , $$
where $r_\nu,  s,t  = 0, \dots, p-1$. We can easily define an 
action of our generators on such states, 
\ba 
    v_\nu |r_1, \dots, r_\nu, \dots, r_{N-1}\rangle _{s,t}
  & = & q^{r_1 + \dots + r_{\nu-1}}
     |r_1, \dots, r_\nu + 1, \dots, r_{N-1}\rangle _{s,t}
      \ \ ,  \nn \\[1mm]
    g_\nu |r_1, \dots, r_\nu, \dots, r_{N-1}\rangle _{s,t}
  & = & q^{r_\nu}
     |r_1, \dots, r_\nu, \dots, r_{N-1}\rangle _{s,t}
      \ \ ,  \nn \\[1mm]
    g_0 |r_1, \dots, r_{N-1} \rangle _{s,t}
  & = & q^{r_1  + \dots + r_{N-1} +  s + t }
     |r_1,  \dots, r_{N-1}\rangle _{s,t}
      \ \ ,  \nn \\[1mm]
    c \ |r_1, \dots, r_{N-1} \rangle _{s,t}
  & = & q^{s}
     |r_1,  \dots, r_{N-1}\rangle _{s,t}
      \ \ .  \nn
\ea
The numerical factor on the right hand side of the  first line is 
an  example of the ``tail of $R$-elements'' discussed above. A similar
term usually appears in the action of monodromies $M$ but is absent here.
It can be seen that the contributions from $(id  \o \vth_N)
(R^{-1})$ and $(id  \o \vth_N)(R)$, which occur in the general 
expression of $D^{IJ}_N(M)$, cancel for $D^{st}_N(c)$, since all 
irreducible representations $\t^r$ of $\Z_q$ are one-dimensional.

{}From the experience with conformal field theories we expect that 
the diagonal representation $\bigoplus D^{\bar I I}_N$ on the 
Hilbert space $\H \equiv \bigoplus_I W^{\bar I I}_N$ is particularly 
relevant. Here we just wish to remark that this representation
can be realized by a very natural construction. Indeed, it was observed
(cf. \cite{AGS}) that algebras such as  $\K_N$ admit a distinguished 
invariant linear functional $\omega: \K_N \rightarrow {\bf C}$, 
$$ \omega( \Je{1}{I_1}_1 \dots \Je{N}{I_N}_N \xi) = \e(\xi) \d_{I_1, 0} 
     \dots \d_{I_N,0} $$
for all $\xi \in \G_n, n \in {\bf Z}\mod N$ and $J^I_n= 
(\t^I \o id)(J_n)$, as usual. When the quantum 
dimensions $d_I$ are positive, this functional is positive and 
hence furnishes -- by the GNS construction -- a Hilbert space 
$\H_\omega$ together with a representation $\pi$ of $\K_N$ on 
$\H_\omega$. If $\vac_\omega$ denotes the GNS vacuum, states 
in $\H_\omega$ are obtained  from 
$$  \Je{1}{I_1}_1 \dots \Je{N}{I_N}_N \vac _\omega \ \ . $$
The formula shows that $\H_\omega$ is isomorphic to the diagonal 
sum $\H = \bigoplus_I W^{\bar I  I}_N \cong \Re^{\o _{N}}$. 
An explicit evaluation of $c^I_{\R,\L}$ on $\H_\omega$ establishes 
an isomorphism of the two spaces  as $\K_N$ modules.     

\begin{prop} The GNS-representation arising from the state 
$\omega: \K_N \rightarrow {\bf C}$ is unitarily equivalent to the 
diagonal representation $\bigoplus_I W^{\bar I I}_N$ of the 
lattice current algebra. 
\end{prop}

The quantum lattice analog of the group-valued local fields
of the WZNW model act in this diagonal representation \cite{LKM,BS}.

\section{Product of Representations}
All continuous current algebras are equipped with a trivial
co-product which
can be written for Fourier modes of currents $j(n)$ as 
\be
\Delta (j(n))=j(n)\otimes 1 + 1\otimes j(n).
\ee
{}From the point of view of  CFT this co-product is not
satisfactory, because it changes the central charge of
representations. In  the framework of CFT, the central charge is 
characteristic of the model, and one must define a new co-product which
preserves it. Such a co-product is provided by the structure of
CFT \cite{Fr}, \cite{MoS}:
\be
\Delta_{CFT}^z(j(n))=j(n)\otimes 1+ 1\otimes \sum_{k\leq n} C^n_k
z^{n-k} j(k).
\ee
Here $C^n_k\equiv n!/k!(n-k)!$ are binomial coefficients.
Observe that the co-product $\Delta_{CFT}$ is not symmetric and
explicitly depends on the parameter $z$.  The aim of this section is
to introduce a lattice counterpart of $\Delta_{CFT}$.

\ssection{Co-product for lattice current algebras}
Because lattice current algebras are labelled by the number of lattice
sites, it is not necessary that both
current algebras on the right hand side of the co-product correspond
to  chains of the same length. We shall define a family of
embeddings
\be
\Lambda_{M,N}:\K_{N+M-1}\rightarrow \K_M \o \K_N
\ee
for any $N$ and $M$. The homomorphisms  $\Lambda_{M,N}$ determine an
action of $\K_{N+M-1}$ on the tensor products $W^{IJ}_M \o W^{KL}_N$
of representation spaces for $\K_M$ and  $\K_N$.

Pictorially, $\Lambda_{M,N}$ corresponds to gluing two closed
chains of length $M$ and $N$ by identifying some site of the first
chain with some site of the second chain.  In this way, the co-product
 $\Lambda_{M,N}$ explicitly depends on the positions of the
identified points. This property is similar to the $z$-dependence of
$\Delta^z_{CFT}$. Below, we always assume that the enumeration
starts from the gluing points, so that this extra parameter does not
show up in our formulas; in the same fashion one can put  {\em e.g.}
 $z=1$ in the continuum theory.

After gluing we cut  the resulting eight-like loop at the middle point and
get one connected chain of  length $N+M-1$. Similarly to
$\Delta^z_{CFT}$, the co-product  $\Lambda_{M,N}$ depends
on the order in which $M$ and $N$ appear. Next, we present the
constructive description of  $\Lambda_{M,N}$.

Let us denote the left currents of $\K_N$ by $J^{\b}_n, n=1,
\dots,N,$
and similarly by $J^{\a}_m, m=1, \dots,M,$ the left currents of
$\K_M$.
$J^{\b}_n$ and $J^{\a}_m$ are regarded as elements in $\G \o \K_M \o
\K_N$
with  the property
$$   \Je{1}{\ \a}_m \Je{2}{\ \b}_n = \Je{2}{\ \b}_n \Je{1}{\ \a}_m
$$
for all  $n,m$. In addition to the left currents, we need $N+M$
commuting  copies of the symmetry  algebra $\G$ to generate
$\K_M  \o \K_N$.

With these notations we can define the announced embedding
$$ \Lambda_{M,N}: \K_{N+M-1} \rightarrow \K_M \o \K_N \ \ . $$
It maps  the generators $J_\rho, \rho =  1,  \dots, N+M-1$,
of $\K_{N+M-1}$ in the following way to generators of 
$\K_M \o \K_N$ :
\be \label{L1}
\Lambda_{M,N}(J_\rho) = \left\{
  \begin{array}{ll}
   J^{\a}_\rho & \mbox{ for } \  \ \rho = 1,  \dots, M-1 \\[1mm]
   J^{\a}_M N^{\a}_-J^{\b}_1 \ \ \  & \mbox{ for } \ \ \rho = M
\\[1mm]
   J^{\b}_{\rho -M +1} & \mbox{ for } \ \ \rho=M+1,
\dots,N+M-2\\[1mm]
   J^{\b}_N (N^{\a}_-)^{-1} & \mbox{ for  } \ \ \rho= N+M-1
  \end{array}  \right.
\ee
where  $N^{\a}_{\pm}=N^{\a}_{0,\pm}$.
For elements  $\xi \in \G$, we define
\be \label{L2}
\Lambda_{N,M}(\iota_\rho(\xi)) = \left\{
  \begin{array}{ll}
   \iota^\a_\rho(\xi)\o e & \mbox{ for } \  \ \rho = 1,  \dots, M-1
\\[1mm]
   e\o \iota^\b_{\rho -M +1}(\xi) & \mbox{ for } \ \ \rho=M,
\dots,N+M-2 \label{Lonxi}\\[1mm]
   (\iota^\a_0\o  \iota^\b_0)(\D(\xi)) & \mbox{ for  } \ \ \rho= N+M-1
  \end{array}\right.
\ee
where  $e$ means the unit element in $\K_M$ or $\K_N$.

\begin{prop} The  map $\Lambda_{M,N}: \K_{N+M-1}  \rightarrow
\K_M \o \K_N$ defined through eqs. (\ref{L1}) and  (\ref{L2}), is
an algebra homomorphism.
\end{prop}
One can prove this proposition by directly verifying the defining
relations for $\K_{N+M-1}$.
In order to justify calling $\Lambda_{M,N}$ a co-product, we need 
some kind of co-associativity.
For lattice current algebras this holds in the following form.

\begin{prop}
The family of homomorphisms $\Lambda_{M,N}$ satisfies the following
property
\be \label{cop}
(id\otimes \Lambda_{N,L})\circ \Lambda_{M, L+N-1}=
(\Lambda_{M,N}\otimes id)\circ \Lambda_{M+N-1,L}.
\ee
\end{prop}

Applying our formulas for $\Lambda_{M,N}$ twice to the generators of
$\K_{L+M+N-2}$ one easily sees that the homomorphisms on the left and on
the right hand sides of (\ref{cop}) coincide.

Let us notice that along with the co-product $\Lambda$ one can introduce
a co-product $\tilde{\Lambda}$ defined by
\be \label{L1'}
\tilde{\Lambda}_{M,N}(J_\rho) = \left\{
  \begin{array}{ll}
   N_+^{\beta}J^{\a}_\rho & \mbox{ for } \  \ \rho = 1 \\[1mm]
   J^{\a}_\rho & \mbox{ for } \  \ \rho = 2, \dots, M-1 \\[1mm]
   J^{\a}_M (N^{\b}_+)^{-1} J^{\b}_1 \ \ \  & \mbox{ for } \ \ \rho = M
\\[1mm]
   J^{\b}_{\rho -M +1} & \mbox{ for } \ \ \rho=M+1,
\dots,N+M-1
  \end{array}  \right.
\ee
and $\tilde {\Lambda}_{M,N}$ acts on elements $\xi \in \G$ according 
to eq. (\ref{Lonxi}) with the $\D$ in the last line being replaced by 
$\D'$. These two co-products $\Lambda$ and $\tilde{\Lambda}$ are 
`intertwined' by the $*$-operation. To make a precise statement we 
introduce an element $K \in \K_M \o \K_N$  by the expression $K = 
\tilde{\Lambda}_{M,N}(\iota_0(\k)) (\iota^\a_0 (\k) \o 
\iota^\b_0(\k))^{-1}$. With this notation we have      
\be
\Lambda_{M,N}(x)^* \ = \ K^{-1}\  \tilde{\Lambda}_{M,N}(x^*) \ K \ \
\ee
for all $x \in \K_{N+M-1}$. This property is similar to (\ref{stop}) 
where the $*$-operation intertwines $\D$ and $\D'$ and, in fact, it 
reduces to the latter on elements $\xi \in \G_n \subset \K_{N+M-1}$. 

\ssection{Special cases}
There are  important special cases of the maps $\Lambda_{M,N}$
which we would like to consider in more detail. First, observe that
the map
\be
 \Lambda_{1,1}: \K_1 \rightarrow \K_1 \o \K_1
\ee
satisfies the co-associativity condition
\be
 (id\otimes \Lambda_{1,1})\circ \Lambda_{1, 1}\ =
\  (\Lambda_{1,1}\otimes id)\circ \Lambda_{1,1}.
\ee
This is a specification of equation (\ref{cop}) for the case of
$L=M=N=1$.  We conclude that the map $\Delta_1=\Lambda_{1,1}$
furnishes a co-product for the algebra  $\K_1$. It has been 
recently shown \cite{FN} that as a Hopf algebra $\K_1$ is 
isomorphic to the Drinfeld double of $\G$. In particular, this
implies existence of an $R$-matrix for the Hopf-algebra $\K_1$. 

As we know (see Section 3), the algebra $\K_1$
is generated by the elements $\xi \in \G$ and by the
universal element $M\in \G \o \K_1$. Irreducible representations
of $\K_1$ are labelled by the pairs $(I, J)$ of  irreducible
representations of $\G$. This implies that as an algebra $\K_1$ is
isomorphic to $\G \o \G$ (see also \cite{ReSTS} where an isomorphism 
of quasi-triangular Hopf algebras is described). 

There is another interesting special choice of chain lengths  $N$ and
$M$:
\be
 \Lambda_{1,N}: \K_N \rightarrow \K_1 \o \K_N\ \ .
\ee
The co-associativity condition (\ref{cop}) adapted to this case reads
\be
 (id\otimes \Lambda_{1,N})\circ \Lambda_{1, N}\ =\ 
 (\Lambda_{1,1}\otimes id)\circ \Lambda_{1,N}\ \ .
\ee
This shows that  $\Lambda_{1,N}$ also provides a co-action of the Hopf
algebra $\K_1$ on $\K_N$. Such a structure has been noticed already
in \cite{LKM}.  We shall see that it permits us to  establish a 
one-to-one correspondence between representations of $\K_N$ and $\K_1$.

Our last remark on the properties of $\Lambda_{M,N}$ concerns the
inductive limit $\K_{\infty}$.
Using the block-spin embeddings $\K_N\rightarrow K_{N+1}$, one can
construct a commutative diagram:

\be \label{d2}
\begin{array}{cccc}
\Lambda_{M,N}: & \K_{N+M-1} & \rightarrow & \K_M \o \K_N \\
                                &\downarrow &                  &
\downarrow \\
\Lambda_{M,N+1}: & \K_{N+M} & \rightarrow & \K_{M} \o \K_{N+1}
\end{array}
\ee
Commutativity (\ref{d2}) ensures that the sequence
of homomorphisms $\Lambda_{M,N}$ defines  homomorphisms
\be
\Lambda_{M,\infty}:\K_{\infty}\rightarrow
\K_M\o \K_{\infty}.
\ee
Then equation (\ref{cop}) implies the co-associativity for
$\Lambda_{M,\infty}$:
\be
(id\otimes \Lambda_{N,\infty})\circ \Lambda_{M\infty}=
 (\Lambda_{M,N}\otimes id)\circ \Lambda_{N+M-1,\infty}.
\ee
Thus, $\Lambda_{M,\infty}$ provides a co-module structure
for $\K_{\infty}$ with respect to the family $\K_M$.
In particular, $\K_{\infty}$ is a co-module over $\K_1$.

\ssection{Implications for  representation theory}
The co-product  $\Lambda_{M,N}$ yields a notion of tensor
product for representations of the algebras $\K_M$ and $\K_N$.

\begin{defn}
The representation $D$  of the algebra $\K_{N+M-1}$ is called a
tensor product of the representations $D_M$ of $\K_M$ and $D_N$ of
$\K_N$ if it acts on the tensor product of the corresponding vector
spaces $W_M\o W_N$ according to the following formula:
\be
D(x) \ = \ (D_M \o D_N) \Lambda_{M,N} (x)\ 
\ee
for all elements $x$ of the lattice current algebra $\K_{N+M-1}$.  
The resulting representation will be denoted by $D_M \bo D_N$. 
\end{defn}
We would like to analyse the structure of this new tensor product.
We denote by $\E$ the trivial representation of the symmetry Hopf
algebra.

\begin{prop}
For any $M$ and $N$ and for arbitrary labels $I$ and $J$ of the
representations of the symmetry algebra the following representations
of $\K_{M+N-1}$ are isomorphic:
\be \label{iso}
D^{IJ}_{M+N-1}\ \simeq \ D^{IJ}_M \bo D^{\E \E}_N
\ \simeq \ D^{\E \E}_N \bo D^{IJ}_M\ \ .
\ee
In this sense, tensoring with the vacuum representation $D^{\E \E}_N $
is trivial.
\end{prop}

To prove this proposition one first checks that on the   spaces
$W^{IJ}_M\o W^{\E \E}_N$ and $W^{\E \E}_N\o W^{IJ}_M$ the central
elements of $\K_{M+N-1}$ have eigenvalues corresponding to the
representation $D^{IJ}_{M+N-1}$. Then one checks that the
dimensions of all these spaces coincide with the dimension of
$W^{IJ}_{M+N-1}$. This completes the proof.

Observe that the tensor product of representations that we have
introduced, relates representations of different algebras. For
instance, if we take $N=M$, we obtain a representation of the
algebra $\K_{2N-1}$ on the tensor product of representation
spaces of $\K_N$. The idea now is to embedd the algebra $\K_N$
into $\K_{2N-1}$ with the help of the block spin maps $\cc_M$,
$$ \cc_{2N-1,N}  \ \ :=\ \  \cc_{2N-2} \circ \dots \circ
    \cc_{N+1} \circ \cc_N : \ \K_N \rightarrow \K_{2N-1}\ \ .
$$
In this way we may represent the algebra $\K_N$ on tensor
products of its own representation spaces. The resulting
representation certainly has a huge commutant and is not
appropriate to describe the representation theory of $\K_N$.
We shall define a certain projection operator $\P^\a_N \in
\K_N \o \K_N$ that projects to more interesting
subrepresentation.

The construction of $\P^\a_N$ proceeds as follows. Notice that
the lattice current algebra $\K_N$ contains $N-1$ local projectors
$p_i \in \G_n \subset \K_N$ where $n$ runs from $1$ to $N-1$.
They are uniquely determined by the property $\iota_i(\xi) p_i
= \e(\xi) p_i$ for all $\xi \in \G$. An explicit formula for
$p_i$ in terms of the objects $N_i$ can be obtained along the
lines of Subsection 2.5. When $N_i$ replaces the monodromy $M$
then we obtain $p_i$ instead of $\chi^0$. These projectors
$p_i$ commute with each other, i.e., $p_i p_j = p_j p_i$ for
all $i,j$ so that their product defines again a projector
$\P_N$,
$$     \P_N \ := \ \prod_{i=1}^{N-1} \ p_i \ \in \K_N \ \ . $$
{}From the defining relations of $\K_N$ it is fairly obvious that
$\P_N$ commutes with $N_0$ and the monodromy $M$. In the
following, $\P^\a_N$ will denote the projector $\P_N \o e \in
\K_N \o \K_N$ and similarly $\P^\b_N = e \o \P_N \in \K_N \o
\K_N$.

With these objects at hand, we are now able to define a new
co-product $\D_N$ for the algebra $K_N$,
\be \label{coprodN}
    \D_N(x) \ : =\  \P^a_N \ \Lambda_{N,N}\, (\cc_{2N-1,N}(x)) \ \
    \mbox{for all} \ \ x \in \K_N \ \ .
\ee
It is easy to check that $\D_N$ defines a homomorphism because
$\cc_{2N-1,N}$ and $\Lambda_{N,N}$ are homomorphisms and
the projector $\P^\a_N$ commutes with the image of $\Lambda_{N,N}
\circ \cc_{2N-1,N} : \K_N \rightarrow \K_N \o \K_N$. Let us also
mention that $\D_N(\P_N \ x ) = \P^\b_N \D_N( x )$ so that the
co-associativity of $\D_N$ follows from that of $\Lambda_{N,M}$,
$$    (\D_N \o id) \D_N (x) \ = \ (id \o \D_N) \D_N (x)
     \ \ \mbox{ for all } \ \ \xi \in \K_N \ \ . $$
In deviation from the standard properties of co-products, $\D_N$
is not unit preserving, i.e., $ \D_N(e)  \neq e \o e \in \K_N
\o \K_N$ ($e \in \K_N$ denotes the unit element), and there is
no one-dimensional trivial representation of $\K_N$. The role
of the co-unit is actually played by the vacuum representation
$D^{\E \E}_N$ of $\K_N$. Such algebraic properties are characteristic
for {\em weak Hopf-algebras} \cite{WHA} and the closely
related  {\em  weak quasi-Hopf algebras} of \cite{MS}.

Notice, that the co-product $\D_N$ is compatible with the block
spin operation:
$(\c_N \otimes \c_N) \D_N = \D_{N+1} \c_N$. This property is
ensured by the fact that two block spin operations $\c_N$ and $\cc_N$
commute with each other (see Section 2). Thus, one can define 
an operation $\D_{\infty}: \K_{\infty} \rightarrow \K_{\infty} \otimes
\K_{\infty}$ which provides a co-product for the inductive limit of
lattice current algebras.

We would finally like to compare the representation category of
$\K_N$ with that of $\K_1$. To this end notice that the formula 
\be   \label{pre}
D^{IJ}_{N}\ \simeq \ D^{IJ}_1 \bo D^{\E \E}_N\ \  
\ee
provides a one-to-one correspondence between representations of
$\K_1$ and $\K_N$ for arbitrary $N$. In fact, this implies the
same kind of correspondence for representations of $\K_1$ and
$\K_{\infty}$. Because $\K_N$ is semisimple, for all $N$, the
isomorphism (\ref{pre}) induces a map
\begin{equation}
\Fu_N(D_1)\ =\ D_N
\end{equation}
which assigns to each representation $D_1$ of the algebra $\K_1$
a representation $D_N$ of the algebra $\K_N$.

To describe the properties of the map $\Fu$, it is convenient to
use the language of the theory of categories (see {\em e.g.}
\cite{Ker}). It is clear that the map $\Fu_N$ is invertible and
that it defines a co-variant tensor functor mapping the category 
of representations of the Hopf-algebra $\K_1$ into the category 
of representations of the lattice current algebras $\K_N$. 
Actually, on the image of $P^\a_N$, the tensor product of 
representations of $\K_N$ defined through $\D_N$ is 
isomorphic to the representation obtained with the help
of the co-action $\Lambda_{1,N}$. Since tensor operators for 
the latter may be trivially identified with tensor operators 
of the quasitriangular Hopf-algebra $\K_1$, the functor $F_N$ 
provides an equivalence of braided tensor categories. In the 
limit $N \to \infty$ we arrive at the following conclusion.

\begin{theo} 
The functor $\Fu$ establishes an isomorphism between representations
the Hopf algebra $\K_1$ and the lattice current algebra $\K_{\infty}$
which is compatible with a co-products of $\K_1$ and $\K_\infty$ and
establishes an equivalence of braided monoidal categories. 
\end{theo}

We can view this fact as the lattice analogue of a theorem in
\cite{KaLu}, \cite{Fin} on the equivalence of tensor categories
corresponding to quantum groups and current algebras. Here the
algebra $\K_{\infty}$ replaces the current algebra, and $\K_1=\G \o
\G$ is the direct product of two quantum groups corresponding to two
chiral sectors. From this point of view, finding an exact relationship
between lattice and continuum current algebras emerges as a
challenging problem.

\vskip 0.2cm \noindent
{\bf Acknowledgements:} We would like to thank A. Connes and
K. Gawedzki, participants and lectures of the 95' Summer school on
Theoretical Physics at Les Houches  for an inspiring atmosphere.
V.S. would also like to thank T. Miwa and I. Ojima for their 
hospitality at RIMS. We are grateful to A. Bytsko and F. Nill 
for their remarks and criticism.


\begin{thebibliography}{99}
\bibitem{Hopf} E. Abe, {\em Hopf algebras}, Cambridge University
 Press, Cambridge, 1980; M.E. Sweedler, {\em Hopf algebras}, Benjamin
 1962
\bibitem{LKM} A. Yu. Alekseev, L.D. Faddeev, M.A.
 Semenov-Tian-Shansky
 {\em Hidden Quantum groups inside Kac-Moody algebras},
 Commun. Math. Phys. {\bf149}, no.2 (1992) p.335\\
 L. D. Faddeev, {\em Quantum symmetry in conformal field theory by
 Hamiltonian methods}, in {\em New symmetry principles in quantum
 field theory}, Proceedings Cargese 91, Plenum Press
\bibitem{A-F} A. Yu. Alekseev, L.D. Faddeev, {\em An evolution and
 dynamics for the $q$-deformed quantum top},  Helsinki-Uppsala preprint,
 December 1994, hep-th/9406196
\bibitem{AGS} A. Yu. Alekseev, H. Grosse, V. Schomerus, {\em
 Combinatorial quantization of the Hamiltonian Chern Simons theory
 I}, hep-th 9403066, Commun. Math. Phys. {\bf 172} (1995) 317-358;\\
 A. Yu. Alekseev, H. Grosse, V. Schomerus, {\em
 Combinatorial quantization of the Hamiltonian Chern Simons theory
 II},hep-th 9408097, Commun. Math. Phys. {\bf 174} (1995) 561-604
\bibitem{AlRe} A.Yu. Alekseev, A. Recknagel, {\em The embedding
 structure and the shift operator of the $U(1)$ lattice current
 algebra}, hep-th/9503107, Lett. Math. Phys. {\bf 37} (1996) 15-27
\bibitem{AlSc} A.Yu. Alekseev, V. Schomerus, {\em Representation
 theory of Chern Simons Observables}, q-alg/9503016, Duke Math. 
 Journal Vol. {\bf 85}, No. 2 (1996) 447 
\bibitem{WHA} G. B\"ohm, K. Szlachanyi, {\em A coassociative 
 $C^*$-quantum group with non-integral dimension}, q-alg/9509008, 
 Lett. Math. Phys. (to appear)
\bibitem{BS} A.Bytsko, V.Schomerus, 
{\em Vertex operators - from a toy model to lattice algebras},
q-alg/9611010
\bibitem{DHR} S. Doplicher, R. Haag, J.E. Roberts, {\em Fields,
 observables and gauge transformations I,II}, Commun. Math. Phys.
 {\bf 13} (1969) 1 and {\bf 15} (1969) 173
\bibitem{DHR2}  S. Doplicher, R. Haag, J.E. Roberts, {\em
 Local observables and particle statistics I,II}, Commun. Math. Phys.
 {\bf 23} (1971) 199 and {\bf 35} (1974) 49   
\bibitem{FRT} L.D. Faddeev, N.Yu. Reshetikhin, L.A.
 Takhtajan, {\em Quantization of Lie Groups and Lie Algebras},
 Algebra and Analysis {\bf 1} (1989) 178-201, 1 and Leningrad Math. J.
 Vol. 1 (1990) 193-225
\bibitem{FaVo} L.D. Faddeev, A.Yu. Volkov, {\em Abelian current algebra
 and the Virasoro  algebra on the lattice}, Phys. Lett. {\bf
 B315} (1993) 311
\bibitem{FaGa} F. Falceto, K. Gawedzki, {\em Lattice
 Wess-Zumino-Witten
 model and quantum groups}, J. Geom. Phys. {\bf 11} (1993) 251
\bibitem{Fin} M. Finkelberg, {\em Fusion Categories},
 Harvard University thesis, May 1993
\bibitem{Fr} J. Fr\"{o}hlich, in "Recent Developments in QFT", Proc.
of 24th Intl.
Conference on High Energy Physics, Munich 1988, Kotthaus \& Kuehn
(eds.)
Springer 1989.
\bibitem{KaLu} D. Kazhdan, G. Lusztig, {\em Tensor structures
 arising from affine Lie algebras I/II}, J. Am. Math. Soc.
 {\bf vol. 6} 4 (1993) 905
\bibitem{Ker} T. Kerler, {\em Nontannakian categories in quantum 
field theory},
Cargese 91, Proceedings, New symmetry principles in quantum
field theory, 492
\bibitem{MS} G. Mack, V. Schomerus, {\em Quasi Hopf quantum
 symmetry in quantum theory}, Nucl. Phys. {\bf B370} (1992),185
\bibitem{Maj} S. Majid, {\em Algebras and Hopf algebras in braided 
 categories}, in {\em Advances in Hopf Algebras}, Lecture Notes in 
 Pure and Appl. Math. {\bf 158}, Dekker, New York, 1994, 55 
\bibitem{MoS} G. Moore, N. Seiberg, {\em Classical  and quantum conformal
 field theory}, Commun. Ma th. Phys. {\bf 123} (1989) 177
\bibitem{FN} F. Nill, {\em On the Structure of Monodromy Algebras and 
 Drinfeld Doubles} q-alg/9609020
\bibitem{NS} F. Nill, K. Szlachanyi, {\em Quantum chains of Hopf algebras
 with quantum double cosymmetry}, preprint hep-th/9509100, Commun. 
 Math. Phys., to appear 
\bibitem{ReSTS} N. Reshetikhin, M. Semenov-Tian-Shansky, {\em
 } Lett. Math. Phys. {\bf 19} (1990) 133; \\
N.Reshetikhin, M.Semenov-Tian-Shansky, {\em Quantum R-matrices
and factorization problems}, J. Geom. Phys. {\bf 5} (1988) 533
\bibitem{ReTu} N. Reshetikhin, V. G. Turaev, {\em Ribbon
 graphs and their invariants derived from quantum groups}, Commun.
 Math. Phys. {\bf 127}, (1990),1; \\
 N. Reshetikhin, V. G. Turaev, {\em Invariants of
 3-manifolds via link polynomials and quantum groups}, Invent.
 Math. {\bf 103} (1991) 547
\bibitem{Sch} V. Schomerus, {\em Construction of field algebras
with quantum symmetry from local observables}, Commun. Math.
Phys. {\bf 169} (1995) 193
\bibitem{Sklyanin-Volkov} E.Sklyanin, A.Volkov, unpublished.
\bibitem{SzVe} K. Szlachanyi, P. Vecsernyes, {\em Quantum symmetry
 and braid group statistics in G-spin models}, Commun. Math. Phys.
 {\bf 156} (1993) 127
\bibitem{Vec} P. Vecsernyes, {\em On the quantum symmetry of the
chiral Ising model}, Nucl. Phys. {\bf B415} (1994) 557
\bibitem{Ver} E. Verlinde, {\em Fusion rules and modular
 transformations in 2D conformal field theory}, Nucl. Phys.
 {\bf B 300} (1988) 360
\bibitem{WZW} J. Wess, B. Zumino, {\em  } Phys. Rev. {\bf B37}
 (1971) 95; E. Witten, {\em } Commun. Math. Phys. {\bf 92} (1984)
 455; S. Novikov, {\em } Uspehi Math. Sci. {\bf 37} (1982) 3
 \end{thebibliography}
\end{document}